# Stochastic Restricted Biased Estimators in misspecified regression model with incomplete prior information

**Manickavasagar Kayanan[1,2] and Pushpakanthie Wijekoon[3]**

[1]Deparment of Physical Science, Vavuniya Campus of the University of Jaffna, Vavuniya, Sri Lanka,
[2] Postgraduate Institute of Science, University of Peradeniya, Peradeniya, Sri Lanka,
[3]Department of Statistics and Computer Science, University of Peradeniya, Peradeniya, Sri Lanka.
*email: mgayanan@vau.jfn.ac.lk, pushpaw@pdn.ac.lk*

In this article, the analysis of misspecification was extended to the recently introduced stochastic restricted biased estimators when multicollinearity exists among the explanatory variables. The Stochastic Restricted Ridge Estimator (SRRE), Stochastic Restricted Almost Unbiased Ridge Estimator (SRAURE), Stochastic Restricted Liu Estimator (SRLE), Stochastic Restricted Almost Unbiased Liu Estimator (SRAULE), Stochastic Restricted Principal Component Regression Estimator (SRPCR), Stochastic Restricted r-k class estimator (SRrk) and Stochastic Restricted r-d class estimator (SRrd) were examined in the misspecified regression model due to missing relevant explanatory variables when incomplete prior information of the regression coefficients is available. Further, the superiority conditions between estimators and their respective predictors were obtained in the mean square error matrix (MSEM) sense. Finally, a numerical example and a Monte Carlo simulation study were used to illustrate the theoretical findings.

## 1 Introduction

Misspecification due to left out relevant explanatory variables is very often when considering the linear regression model, which causes these variables to become a part of the error term. Consequently, the expected value of error term of the model will not be zero. Also, the omitted variables may be correlated with the variables in the model. Therefore, one or more assumptions of the linear regression model will be violated when the model is misspecified, and hence the estimators become biased and inconsistent. Further, it is well-known that the ordinary least squares estimator (OLSE) may not be very reliable if multicollinearity exists in the linear regression model. As a remedial measure to solve multicollinearity problem, biased estimators based on the sample model $y = X\beta + \varepsilon$ with prior information which can be exact or stochastic restrictions have received much attention in the statistical literature. The intention of this work is to examine the performance of the recently introduced stochastic restricted biased estimators in the misspecified regression model with incomplete prior knowledge about regression coefficients when there exists multicollinearity among explanatory variables.

When we consider the biased estimation in misspecified regression model without any restrictions on regression parameters, Sarkar (1989) discussed the consequences of exclusion of some important explanatory variables from a linear regression model when multicollinearity exists. Şiray (2015) and Wu (2016) examined the efficiency of the r-d class estimator and r-k class estimator over some existing estimators, respectively in the misspecified regression model. Chandra and Tyagi (2017) studied the effect of misspecification due to the omission of relevant variables on the dominance of the r-(k,d) class estimator. Recently, Kayanan and Wijekoon (2017) examined the performance of existing biased estimators and the respective predictors based on the sample information in a misspecified linear regression model without considering any prior information about regression coefficients.

It is recognized that the mixed regression estimator (MRE) introduced by Theil and Goldberger (1961) outperform ordinary least squares estimator (OLSE) when the regression model is correctly specified. The biased estimation with stochastic linear restrictions in the misspecified regression model due to inclusion of an irrelevant variable with the incorrectly specified prior information was discussed by Teräsvirta (1980). Later Mittelhmmer (1981), Ohtani and Honda (1984), Kadiyala (1986) and Trenkler and Wijekoon (1989) discussed the efficiency of MRE under misspecified regression model due to exclusion of a relevant variable with correctly specified prior information. Further, the superiority of MRE over the OLSE under the misspecified regression model with incorrectly specified sample and prior information was discussed by Wijekoon and Trenkler (1989). Hubert and Wijekoon (2004) have considered the improvement of Liu estimator (LE) under a misspecified regression model with stochastic restrictions, and introduced the Stochastic Restricted Liu Estimator (SRLE).

In this paper, the performance of the recently introduced stochastic restricted estimators namely the Stochastic Restricted Ridge Estimator (SRRE) proposed by Li and Yang (2010), Stochastic Restricted Almost Unbiased Ridge Estimator (SRAURE) and Stochastic Restricted Almost Unbiased Liu Estimator (SRAULE) proposed by Wu and Yang (2014), Stochastic Restricted Principal Component Regression Estimator (SRPCR) proposed by He and Wu (2014), Stochastic Restricted r-k class estimator (SRrk) and Stochastic Restricted r-d class estimator (SRrd) proposed by Jibo Wu (2014) were examined in the misspecified regression model when multicollinearity exists among explanatory variables. Further, a generalized form to represent these estimators is also proposed.

The rest of this article is organized as follows. The model specification and the estimators are written in section 2. In section 3, the Mean Square Error Matrix (MSEM) comparison between two estimators and respective predictors are considered. In section 4, a numerical example and a Monte Carlo simulation study are given to illustrate the theoretical results in Scalar Mean Square Error (SMSE) criterion. Finally, some concluding remarks are mentioned in section 5. The references and appendixes are given at the end of the paper.

## 2 Model specification and the estimators

Assume that the true regression model is given by
$$y = X_1\beta_1 + X_2\beta_2 + \varepsilon = X_1\beta_1 + \delta + \varepsilon \tag{2.1}$$
where $y$ is the $n \times 1$ vector of observations on the dependent variable, $X_1$ and $X_2$ are the $n \times l$ and $n \times p$ matrices of observations on the $m = l + p$ regressors, $\beta_1$ and $\beta_2$ are the $l \times 1$ and $p \times 1$ vectors of unknown coefficients, $\varepsilon$ is the $n \times 1$ vector of disturbances such that $E(\varepsilon) = 0$ and $E(\varepsilon\varepsilon') = \Omega = \sigma^2 I$. Let us say that the researcher misspecifies the regression model by excluding $p$ regressors as
$$y = X_1\beta_1 + u \tag{2.2}$$
Let us also assume that there exists prior information on $\beta_1$ in the form of
$$r = R\beta_1 + g + v \tag{2.3}$$
where $r$ is the $q \times 1$ vector, $R$ is the given $q \times l$ matrix with rank $q$, $g$ is the $q \times 1$ unknown fixed vector, $v$ is the $q \times 1$ vector of disturbances such that $E(v) = 0$, $D(v) = E(vv') = \Psi = \sigma^2 W$, where $W$ is positive definite, and $E(vu') = 0$

By combining sample model (2.2) and prior information (2.3), Thiel and Goldberger (1961) proposed the Mixed Regression Estimator (MRE) as
$$\hat{\beta}_{MRE} = (X_1'\Omega^{-1}X_1 + R'\Psi^{-1}R)^{-1}(X_1'\Omega^{-1}y + R'\Psi^{-1}r) \tag{2.4}$$
$$= (X_1'X_1 + R'W^{-1}R)^{-1}(X_1'y + R'W^{-1}r)$$

To combat multicollinearity, several researchers introduce different types of stochastic restricted estimators in place of MRE. Seven such estimators are SRRE, SRAURE, SRLE, SRALUE, SRPCRE, SRrk class estimator and SRrd class estimator defined below, respectively:
$$\hat{\beta}_{SRRE} = (X_1'X_1 + kI)^{-1}X_1'X_1\hat{\beta}_{MRE} \tag{2.5}$$

$$\hat{\beta}_{SRAURE} = (I - k^2(X_1'X_1 + kI)^{-2})\hat{\beta}_{MRE} \tag{2.6}$$
$$\hat{\beta}_{SRLE} = (X_1'X_1 + I)^{-1}(X_1'X_1 + dI)\hat{\beta}_{MRE} \tag{2.7}$$
$$\hat{\beta}_{SRAULE} = (I - (1-d)^2(X_1'X_1 + I)^{-2})\hat{\beta}_{MRE} \tag{2.8}$$
$$\hat{\beta}_{SRPCR} = T_h T_h' \hat{\beta}_{MRE} \tag{2.9}$$
$$\hat{\beta}_{SRrk} = T_h T_h'(X_1'X_1 + kI)^{-1}X_1'X_1 \hat{\beta}_{MRE} \tag{2.10}$$
$$\hat{\beta}_{SRrd} = T_h T_h'(X_1'X_1 + I)^{-1}(X_1'X_1 + dI)\hat{\beta}_{MRE} \tag{2.11}$$

where $k > 0$, $0 < d < 1$ and $T_h = (t_1, t_2, \ldots, t_h)$ be the first $h$ coloumns of $T = (t_1, t_2, \ldots, t_h, \ldots, t_l)$ which is an orthogonal matrix of the standardized eigenvectors of $X_1'X_1$.

According to Kadiyala (1986), now we apply the simultaneous decomposition to the two symmetric matrices $X_1'X_1$ and $R'\Psi^{-1}R$, as
$$B'X_1'X_1 B = I \quad \text{and} \quad B'R'\Psi^{-1}RB = \Lambda,$$
where $X_1'X_1$ is a positive definite matrix and $R'\Psi^{-1}R$ is a positive semi-definite matrix, $B$ is a $l \times l$ nonsingular matrix, $\Lambda$ is a $l \times l$ diagonal matrix with eigenvalues $\lambda_i > 0$ for $i = 1, 2, \ldots, q$ and $\lambda_i = 0$ for $i = q+1, \ldots, l$.

Let $X_* = X_1 B$, $R_* = RB$, $\gamma = B^{-1}\beta_1$, $X_*'X_* = I$ and $R_*'\Psi^{-1}R_* = \Lambda$, then the models (2.1), (2.3) and (2.3) can be written as

$$y = X_*\gamma + \delta + \varepsilon, \tag{2.12}$$

$$y = X_*\gamma + u, \tag{2.13}$$

$$r = R_*\gamma + g + v. \tag{2.14}$$

According to Wijekoon and Trenkler (1989), the corresponding MRE is given by

$$\hat{\gamma}_{MRE} = (X_*'X_* + R_*'\Psi^{-1}R_*)^{-1}(X_*'y + R_*'W^{-1}r)$$
$$= (I + \sigma^2 \Lambda)^{-1}(X_*'y + R_*'W^{-1}r) \tag{2.15}$$

Hence, the respective expectation vector, bias vector and dispersion matrix are given by

$$E(\hat{\gamma}_{MRE}) = \gamma + (I + \sigma^2 \Lambda)^{-1}(X_*'\delta + R_*'W^{-1}g), \tag{2.16}$$
$$Bias(\hat{\gamma}_{MRE}) = (I + \sigma^2 \Lambda)^{-1}(X_*'\delta + R_*'W^{-1}g), \tag{2.17}$$
$$D(\hat{\gamma}_{MRE}) = \sigma^2(I + \sigma^2 \Lambda)^{-1}. \tag{2.18}$$

In the case of misspecification, now the SRRE, SRAURE, SRLE, SRAULE, SRPCRE, SRrk and SRrd for model (2.12) can be written as

$$\hat{\gamma}_{SRRE} = (X_*'X_* + kI)^{-1}X_*'X_* \hat{\gamma}_{MRE} = (1+k)^{-1}\hat{\gamma}_{MRE} = C_k \hat{\gamma}_{MRE} \tag{2.19}$$
$$\hat{\gamma}_{SRAURE} = (I - k^2(X_*'X_* + kI)^{-2})\hat{\gamma}_{MRE} = (1+k)^{-2}(1+2k)\hat{\gamma}_{MRE}$$
$$= (1+2k)(C_k)^2 \hat{\gamma}_{MRE} = C_k^* \hat{\gamma}_{MRE} \tag{2.20}$$
$$\hat{\gamma}_{SRLE} = (X_*'X_* + I)^{-1}(X_*'X_* + dI)\hat{\gamma}_{MRE} = 2^{-1}(1+d)\hat{\gamma}_{MRE} = C_d \hat{\gamma}_{MRE} \tag{2.21}$$
$$\hat{\gamma}_{SRAULE} = (I - (1-d)^2(X_*'X_* + I)^{-2})\hat{\gamma}_{MRE} = 2^{-2}(1+d)(3-d)\hat{\gamma}_{MRE}$$
$$= 2^{-1}(3-d)C_d \hat{\gamma}_{MRE} = C_d^* \hat{\gamma}_{MRE} \tag{2.22}$$
$$\hat{\gamma}_{SRPCR} = T_h T_h' \hat{\gamma}_{MRE} = C_h \hat{\gamma}_{MRE} \tag{2.23}$$
$$\hat{\gamma}_{SRrk} = (1+k)^{-1} T_h T_h' \hat{\gamma}_{MRE} = C_k C_h \hat{\gamma}_{MRE} = C_{hk} \hat{\gamma}_{MRE} \tag{2.24}$$
$$\hat{\gamma}_{SRrd} = 2^{-1}(1+d) T_h T_h' \hat{\gamma}_{MRE} = C_d C_h \hat{\gamma}_{MRE} = C_{hd} \hat{\gamma}_{MRE} \tag{2.25}$$

respectively, where $C_k = (1+k)^{-1}$, $C_k^* = (1+2k)(C_k)^2$, $C_d = 2^{-1}(1+d)$, $C_d^* = 2^{-1}(3-d)C_d$, $C_h = T_r T_r'$, $C_{hk} = C_k C_h$ and $C_{hd} = C_d C_h$.

It is clear that $C_k$, $C_k^*$, $C_d$ and $C_d^*$ are positive definite, and $C_h$, $C_{hk}$ and $C_{hd}$ are non-negative definite.

Since all these estimators can be written by incorporating $\hat{\gamma}_{MRE}$, now we write a generalized form to represent SRRE, SRAURE, SRLE, SRAULE, SRPCR, SRrk and SRrd as given below:

$$\hat{\gamma}_{(j)} = G_{(j)} \hat{\gamma}_{MRE} \quad (2.26)$$

where $G_{(j)}$ is positive definite matrix if it stands for $C_k$, $C_k^*$, $C_d$ and $C_d^*$, and it is non-negative definite matrix if it stands for $C_h$, $C_{hk}$ and $C_{hd}$.

Now the expectation vector, bias vector, the dispersion matrix and the mean square error matrix can be written as

$$E(\hat{\gamma}_{(j)}) = G_{(j)} E(\hat{\gamma}_{MRE}) = G_{(j)}(\gamma + (I+\sigma^2\Lambda)^{-1}(X_*'\delta + R_*'W^{-1}g)) = G_{(j)}(\gamma + \tau A) \quad (2.27)$$

$$Bias(\hat{\gamma}_{(j)}) = E(\hat{\gamma}_{(j)} - \gamma) = G_{(j)}(\gamma + \tau A) - \gamma = (G_{(j)} - I)\gamma + G_{(j)}\tau A \quad (2.28)$$

$$D(\hat{\gamma}_{(j)}) = G_{(j)} D(\hat{\gamma}_{MRE}) G_{(j)}' = \sigma^2 G_{(j)}(I+\sigma^2\Lambda)^{-1} G_{(j)}' = \sigma^2 G_{(j)} \tau G_{(j)}' \quad (2.29)$$

$$\begin{aligned} MSEM(\hat{\gamma}_{(j)}) &= E(\hat{\gamma}_{(j)} - \gamma)(\hat{\gamma}_{(j)} - \gamma)' \\ &= D(\hat{\gamma}_{(j)}) + Bias(\hat{\gamma}_{(j)}) Bias(\hat{\gamma}_{(j)})' \\ &= \sigma^2 G_{(j)} \tau G_{(j)}' + \left((G_{(j)} - I)\gamma + G_{(j)}\tau A\right)\left((G_{(j)} - I)\gamma + G_{(j)}\tau A\right)' \end{aligned} \quad (2.30)$$

where $\tau = (I+\sigma^2\Lambda)^{-1}$ and $A = (X_*'\delta + R_*'W^{-1}g)$.

Based on 2.27 to 2.30, the respective bias vector, dispersion matrix and MSEM of the MRE, SRRE, SRAURE, SRLE, SRAULE, SRPCR, SRrk and SRrd can easily be obtained, and given in Table B1 in Appendix B.

By using the approach of Kadiyala (1986) and equations (2.3) and (2.4), the generalized prediction function can be defined as follows:

$$y_0 = X_*\gamma + \delta \quad (2.31)$$
$$\hat{y}_{(j)} = X_* \hat{\gamma}_{(j)} \quad (2.32)$$

where $y_0$ is the actual value and $\hat{y}_{(j)}$ is the corresponding predictor.

The MSEM of the generalized predictor is given by

$$\begin{aligned} MSEM(\hat{y}_{(j)}) &= E(\hat{y}_{(j)} - y_0)(\hat{y}_{(j)} - y_0)' \\ &= X_*\left(MSEM(\hat{\gamma}_{(j)})\right)X_*' - X_*\left(Bias(\hat{\gamma}_{(j)})\right)\delta' - \delta\left(Bias(\hat{\gamma}_{(j)})\right)'X_*' + \delta\delta' \end{aligned} \quad (2.33)$$

Note that the predictors based on the MRE, SRRE, SRAURE, SRLE, SRAULE, SRPCR, SRrk and SRrd are denoted by $\hat{y}_{MRE}$, $\hat{y}_{SRRE}$, $\hat{y}_{SRAURE}$, $\hat{y}_{SRLE}$, $\hat{y}_{SRAULE}$, $\hat{y}_{SRPCR}$, $\hat{y}_{SRrk}$ and $\hat{y}_{SRrd}$ respectively.

## 3 Mean Square Error Matrix (MSEM) comparisons

If two generalized biased estimators $\hat{\gamma}_{(i)}$ and $\hat{\gamma}_{(j)}$ are given, the estimator $\hat{\gamma}_{(j)}$ is said to be superior to $\hat{\gamma}_{(i)}$ with respect to MSEM sense if and only if $MSEM(\hat{\gamma}_{(i)}) - MSEM(\hat{\gamma}_{(j)}) \geq 0$. Also, if two generalized predictors $\hat{y}_{(i)}$ and $\hat{y}_{(j)}$ are given, the predictor $\hat{y}_{(j)}$ is said to be superior to $\hat{y}_{(i)}$ with respect to MSEM sense if and only if $MSEM(\hat{y}_{(i)}) - MSEM(\hat{y}_{(j)}) \geq 0$.

Now let $D_{(i,j)} = D(\hat{\gamma}_{(i)}) - D(\hat{\gamma}_{(j)})$, $b_{(i)} = Bias(\hat{\gamma}_{(i)})$, $b_{(j)} = Bias(\hat{\gamma}_{(j)})$ and

$$\Delta_{(i,j)} = MSEM(\hat{\gamma}_{(i)}) - MSEM(\hat{\gamma}_{(j)}) = D_{(i,j)} + b_{(i)}b'_{(i)} - b_{(j)}b'_{(j)}.$$

By applying Lemma A1 (see Appendix A), the following theorem can be stated for the superiority of $\hat{\gamma}_{(j)}$ over $\hat{\gamma}_{(i)}$ with respect to the MSEM criterion.

**Theorem 1:** If $D_{(i,j)}$ is positive definite, then $\hat{\gamma}_{(j)}$ is superior to $\hat{\gamma}_{(i)}$ in MSEM sense when the regression model is misspecified due to excluding relevant variables if and only if
$$b'_{(j)}(D_{(i,j)} + b_{(i)}b'_{(i)})^{-1}b_{(j)} \leq 1.$$

**Proof:** Let $D_{(i,j)}$ is a positive definite matrix. According to Lemma A1 (see Appendix A), $\Delta_{(i,j)}$ is non-negative definite matrix if $b'_{(j)}(D_{(i,j)} + b_{(i)}b'_{(i)})^{-1}b_{(j)} \leq 1$. This completes the proof.

The following theorem can be stated for the superiority of $\hat{y}_{(j)}$ over $\hat{y}_{(i)}$ with respect to the MSEM criterion.

**Theorem 2:** If $A \geq 0$, $\hat{y}_{(j)}$ is superior to $\hat{y}_{(i)}$ in MSEM sense when the regression model is misspecified due to excluding relevant variables if and only $\theta \in \Re(A)$ and $\theta' A^{-1}\theta \leq 1$, where $A = X_*\Delta_{(i,j)}X'_* + X_*(b_{(i)} - b_{(j)})(b_{(i)} - b_{(j)})'X'_* + \delta\delta'$, $\theta = \delta + X_*(b_{(i)} - b_{(j)})$ and $\Re(A)$ stands for column space of $A$ and $A^{-1}$ is an independent choice of g-inverse of $A$.

**Proof:** According to 2.33 we can write $MSEM(\hat{y}_{(i)}) - MSEM(\hat{y}_{(j)})$ as
$$\begin{aligned} MSEM&(\hat{y}_{(i)}) - MSEM(\hat{y}_{(j)}) \\ &= X_*\left(MSEM(\hat{\gamma}_{(i)}) - MSEM(\hat{\gamma}_{(j)})\right)X'_* - X_*\left(Bias(\hat{\gamma}_{(i)}) - Bias(\hat{\gamma}_{(j)})\right)\delta' \\ &\quad - \delta\left(Bias(\hat{\gamma}_{(i)}) - Bias(\hat{\gamma}_{(j)})\right)'X'_* \\ &= X_*\Delta_{(i,j)}X'_* - X_*(b_{(i)} - b_{(j)})\delta' - \delta(b_{(i)} - b_{(j)})'X'_* \end{aligned}$$

After some straight forward calculation, it can be written as
$$MSEM(\hat{y}_{(i)}) - MSEM(\hat{y}_{(j)}) = A - \theta\theta'$$
where $A = X_*\left(\Delta_{(i,j)} + (b_{(i)} - b_{(j)})(b_{(i)} - b_{(j)})'\right)X'_* + \delta\delta'$ and $\theta = \delta + X_*(b_{(i)} - b_{(j)})$.

Due to Lemma A3 (see Appendix A), $MSEM(\hat{y}_{(i)}) - MSEM(\hat{y}_{(j)})$ is non-negative definite matrix if and only if $A \geq 0$, $\theta \in \Re(A)$ and $\theta'A^{-1}\theta \leq 1$, where $\Re(A)$ stands for column space of $A$ and $A^{-1}$ is an independent choice of g-inverse of $A$. This completes the proof.

Based on Theorem 1 and Theorem 2 we can define Corollaries C1-C28, written in the Appendix C, for the superiority conditions between two selected estimators and for the respective predictors by substituting the relevant expressions for $Bias(\hat{\gamma}_{(i)})$, $Bias(\hat{\gamma}_{(j)})$, $D(\hat{\gamma}_{(i)})$ and $D(\hat{\gamma}_{(j)})$ given in Table B1 in Appendix B.

# 4 Illustration of theoretical results
## 4.1 Numerical example

To illustrate the theoretical results, the dataset which gives total National Research and Development Expenditures—as a Percent of Gross National Product by Country: 1972–1986 is considered. The dependent variable $Y$ of this dataset is the percentage spent by the United States, and the four other independent variables are $X_1$, $X_2$, $X_3$ and $X_4$. The variable $X_1$ represents the percent spent by the former Soviet Union, $X_2$ that spent by France, $X_3$ that spent by West Germany, and $X_4$ that spent by the Japan. The data has been analysed by Gruber (1998), Akdeniz and Erol (2013), Li and Yang (2010) and among others. Now we assemble the data as follows:

$$X = \begin{pmatrix} 1.9 & 2.2 & 1.9 & 3.7 \\ 1.8 & 2.2 & 2.0 & 3.8 \\ 1.8 & 2.4 & 2.1 & 3.6 \\ 1.8 & 2.4 & 2.2 & 3.8 \\ 2.0 & 2.5 & 2.3 & 3.8 \\ 2.1 & 2.6 & 2.4 & 3.7 \\ 2.1 & 2.6 & 2.6 & 3.8 \\ 2.2 & 2.6 & 2.6 & 4.0 \\ 2.3 & 2.8 & 2.8 & 3.7 \\ 2.3 & 2.7 & 2.8 & 3.8 \end{pmatrix} \quad y = \begin{bmatrix} 2.3 \\ 2.2 \\ 2.2 \\ 2.3 \\ 2.4 \\ 2.5 \\ 2.6 \\ 2.6 \\ 2.7 \\ 2.7 \end{bmatrix}$$

Note that the eigenvalues of the $X'X$ are 302.96, 0.728, 0.044, 0.035, the condition number is 93, and the variance Inflation Factor (VIF) values are 6.91, 21.58, 29.75, and 1.79. This implies the existence of serious multicollinearity in the data set.

The corresponding OLS estimator of $\beta$ is $\hat{\beta} = (X'X)^{-1}X'y = (0.645, 0.089, 0.143, 0.152)$ and the estimate of $\sigma^2$ is $\hat{\sigma}^2 = 0.00153$. In this example we consider $R = (1, -2, -2, -2)$ and $g = c(1, -1, 2, 0)$. The SMSE values of the estimators are summarized in the Tables B2-B3 in Appendix B.

Table B2 shows the estimated SMSE values of MRE, SRRE, SRAURE, SRLE, SRAULE, SRPCR, SRrk and SRrd for the regression model when $(l, p) = (4, 0)$, $(l, p) = (3, 1)$, and $(l, p) = (2, 2)$ with respect to shrinkage parameters (k/d), where $l$ denotes the number of variable in the model and $p$ denotes the number of misspecified variables. Table B3 shows the estimated SMSE values of the predictor of MRE, SRRE, SRAURE, SRLE, SRAULE, SRPCR, SRrk and SRrd for the regression model when $(l, p) = (4, 0)$, $(l, p) = (3, 1)$, and $(l, p) = (2, 2)$ for some selected shrinkage parameters (k/d).

Note that when $(l, p) = (4, 0)$ the model is correctly specified, when $(l, p) = (3, 1)$ one variable is omitted from the model and when $(l, p) = (2, 2)$ two variables are omitted from the model. For simplicity we choose shrinkage parameter values k and d in the range (0, 1).

From Table B2, we can observe that the MRE is superior to the other estimators when $(l, p) = (4, 0)$, and SRAULE, SRRE, SRLE and SRAURE are outperformed the other estimators for $(k/d) < 0.2$, $0.2 \leq (k/d) < 0.5$, $0.5 \leq (k/d) < 0.7$ and $(k/d) \geq 0.7$, respectively, when $(l, p) = (3, 1)$. Similarly, SRLE and SRRE are superior to the other estimators for $(k/d) < 0.5$ and $(k/d) \geq 0.5$, respectively, when $(l, p) = (2, 2)$.

From Table B3, we further observe that predictors based on SRLE and SRRE are outperformed the other predictors for $(k/d) < 0.5$ and $(k/d) \geq 0.5$, respectively, when $(l, p) = (4, 0)$ and $(l, p) = (3, 1)$, and predictors based on SRrd and SRrk are superior to the other predictors for $(k/d) < 0.5$ and $(k/d) \geq 0.5$, respectively, when $(l, p) = (2, 2)$.

## 4.2 Simulation

For further clarification, a Monte Carlo simulation study is done at different levels of misspecification using R 3.2.5. Following McDonald and Galarneau (1975), we can generate the explanatory variables as follows:

$$x_{ij} = (1 - \rho^2)^{1/2} z_{ij} + \rho z_{i,m} \;; i = 1, 2, \ldots, n. \; j = 1, 2, \ldots, m.$$

where $z_{ij}$ is an independent standard normal pseudo random number, and $\rho$ is specified so that the theoretical correlation between any two explanatory variables is given by $\rho^2$. A dependent variable is generated by using the following equation

$$y_i = \beta_1 x_{i1} + \beta_2 x_{i2} + \beta_3 x_{i3} + \beta_4 x_{i4} + \beta_5 x_{i5} + \varepsilon_i \;; i = 1, 2, \ldots, n.$$

where $\varepsilon_i$ is a normal pseudo random number with mean zero and variance one. Also, we select $\beta = (\beta_1, \beta_2, \beta_3, \beta_4, \beta_5)$ as the normalized eigenvector corresponding to the largest eigenvalue of $X'X$ for which $\beta'\beta = 1$. Further we choose $R = (1, 1, 1, 1, 1)$ and $g = (1, -2, 0, 3, 1)$.

Then the following set up is considered to investigate the effects of different degrees of multicollinearity on the estimators:

- $\rho = 0.9$, condition number =9.49 and VIF = (5.99, 5.88, 5.94, 5.96, 20.47)
- $\rho = 0.99$, condition number =34.77 and VIF = (57.66, 56.50, 57.26, 57.31, 225.06)
- $\rho = 0.999$, condition number =115.66 and VIF = (574.3, 562.8, 570.7, 570.8, 2271.4)

Three different sets of observations are considered by selecting $(l, p) = (5, 0)$, $(l, p) = (4, 1)$ and $(l, p) = (3, 2)$ when $n = 50$, where $l$ denotes the number of variable in the model and $p$ denotes the number of misspecified variables. Note that when $(l, p) = (5, 0)$ the model is correctly specified, when $(l, p) = (4, 1)$ one variable is omitted from the model and when $(l, p) = (3, 2)$ two variables are omitted from the model. For simplicity, we select values $k$ and $d$ in the range $(0,1)$.

The simulation is repeated 2000 times by generating new pseudo random numbers and the simulated SMSE values of the estimators and predictors are obtained using the following equations:

$$SMSE(\hat{\gamma}_{(j)}) = \frac{1}{2000} \sum_{r=1}^{2000} tr\left(MSEM(\hat{\gamma}_{(j)r})\right) \text{ and}$$

$$SMSE(\hat{y}_{(j)}) = \frac{1}{2000} \sum_{r=1}^{2000} tr\left(MSEM(\hat{y}_{(j)r})\right) \text{ respectively.}$$

The simulation results are summarized in Tables B4-B9 in Appendix B.

Table B4, Table B5, and Table B6 in Appendix show the estimated SMSE values of the estimators for the regression model when $(l, p) = (5, 0)$, $(l, p) = (4, 1)$ and $(l, p) = (3, 2)$, and $\rho = 0.9$, $\rho = 0.99$ and $\rho = 0.999$ for the selected values of shrinkage parameters (k/d), respectively. Table B7, Table B8, and Table B9 in Appendix show the corresponding estimated SMSE values of the predictors for the above regression models, respectively.

From Table B4, we can observe that MRE and SRAULE are outperformed the other estimators for $(k/d) < 0.8$ and $(k/d) \geq 0.8$, respectively, when $(l, p) = (5, 0)$ and $(l, p) = (4, 1)$. Further, SRLE and SRRE are superior to the other estimators for $(k/d) < 0.5$ and $(k/d) \geq 0.5$, respectively, when $(l, p) = (3, 2)$ under $\rho = 0.9$.

From Table B5, we can observe that SRAULE, MRE and SRAURE are outperformed the other estimators for $(k/d) < 0.3$, $0.3 \leq (k/d) < 0.7$ and $(k/d) \geq 0.7$, respectively, when $(l, p) = (5, 0)$. Similarly, SRAULE, SRRE, SRLE and SRAURE are superior to the other estimators when $(k/d) < 0.2$, $0.2 \leq (k/d) < 0.5$, $0.5 \leq (k/d) < 0.7$ and $(k/d) \geq 0.7$, respectively, when $(l, p) = (4, 1)$, and both SRLE and SRRE are outperformed the other estimators for $(k/d) < 0.5$ and $(k/d) \geq 0.5$, respectively, when $(l, p) = (3, 2)$ and $\rho = 0.99$.

The results in Table B6 indicate that MRE is superior to the other estimators when $(l, p) = (5, 0)$, and SRAULE, SRRE, SRLE and SRAURE are outperformed the other estimators for $(k/d) < 0.2$, $0.2 \leq (k/d) < 0.5$, $0.5 \leq (k/d) < 0.7$ and $(k/d) \geq 0.7$, respectively, when $(l, p) = (4, 1)$. Further, SRLE and SRRE are outperformed the other estimators for $(k/d) < 0.5$ and $(k/d) \geq 0.5$, respectively, when $(l, p) = (3, 2)$ and $\rho = 0.999$.

From Tables B7-B9, we further observe that the predictors based on SRrd and SRrk are always outperformed the other predictors for $(k/d) < 0.5$ and $(k/d) \geq 0.5$, respectively, when $(l, p) = (5, 0)$, $(l, p) = (4, 1)$ and $(l, p) = (3, 2)$.

The SMSE values of the selected estimators are plotted with different ρ values to demonstrate the results graphically when $(l, p) = (3, 2)$. Figures 1-3 show the graphical illustration of the performance of estimators in the misspecified regression model $((l, p) = (3, 2))$ when $\rho = 0.9$, $\rho = 0.99$ and $\rho = 0.999$, respectively. Similarly, Figures 4-6 present the graphical illustration of the performance of predictors in the misspecified regression model $((l, p) = (3, 2))$ when $\rho = 0.9$, $\rho = 0.99$ and $\rho = 0.999$, respectively.

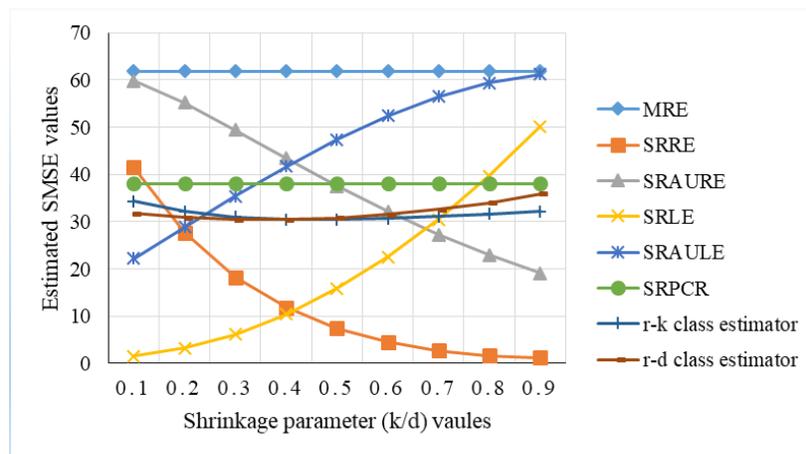

Figure 1. SMSE values of the estimators in the misspecified regression model $((l, p) = (3, 2))$ when n = 50 and $\rho = 0.9$

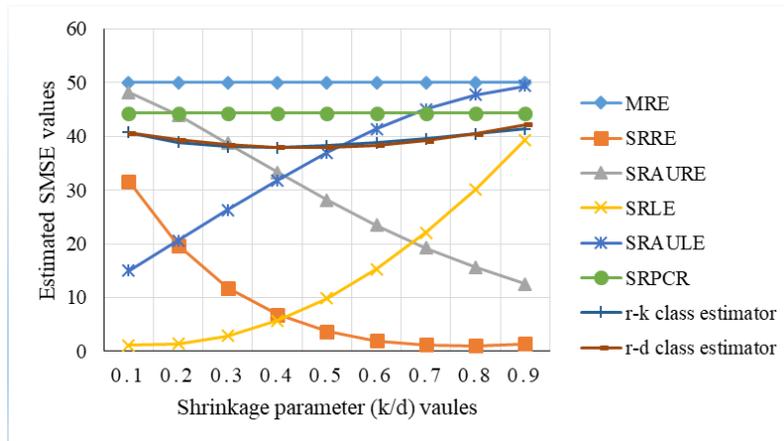

Figure 2. SMSE values of the estimators in the misspecified regression model $((l, p) = (3, 2))$ when n = 50 and $\rho = 0.99$

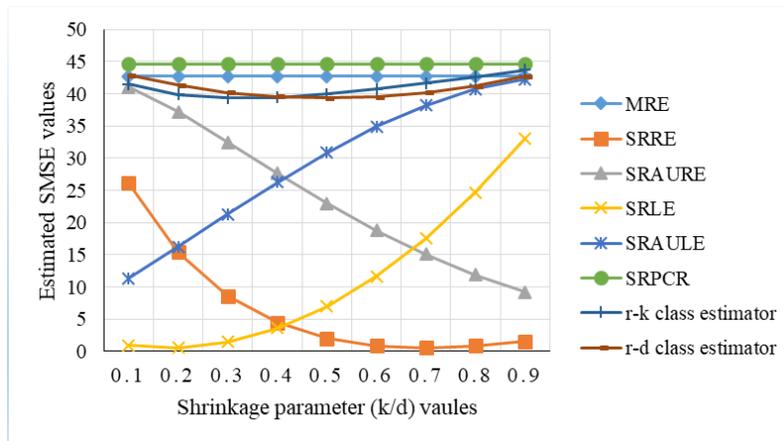

Figure 3. SMSE values of the estimators in the misspecified regression model $((l, p) = (3, 2))$ when n = 50 and $\rho = 0.999$

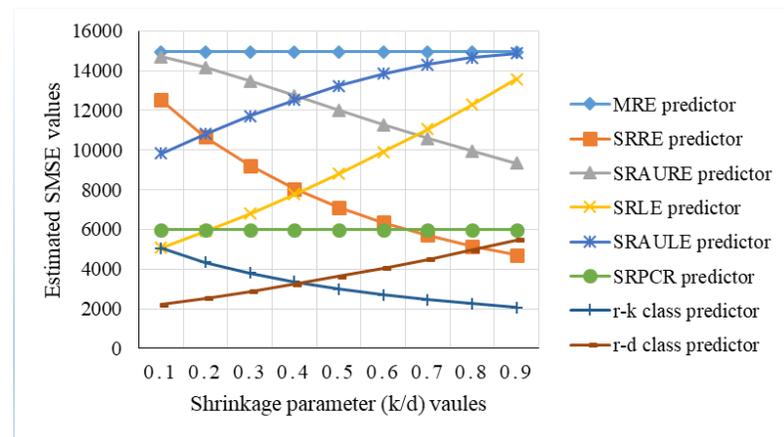

Figure 4. SMSE values of the predictors in the misspecified regression model $((l, p) = (3, 2))$ when n = 50 and $\rho = 0.9$

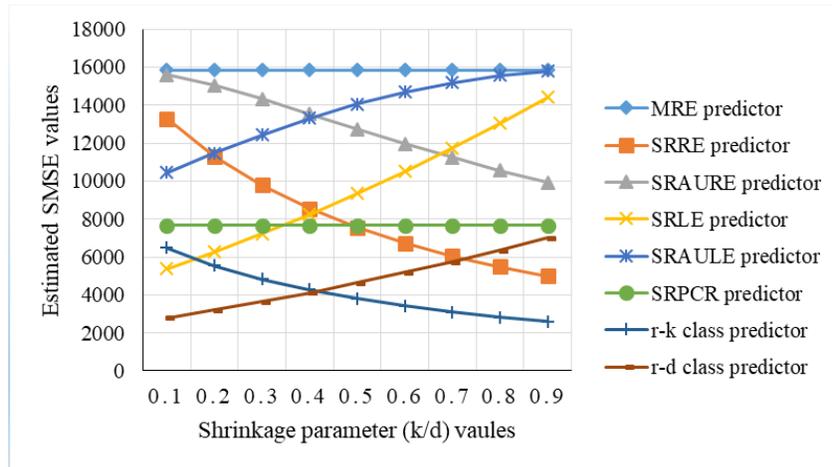

Figure 5. SMSE values of the predictors in the misspecified regression model $((l, p) = (3, 2))$ when n = 50 and $\rho = 0.99$

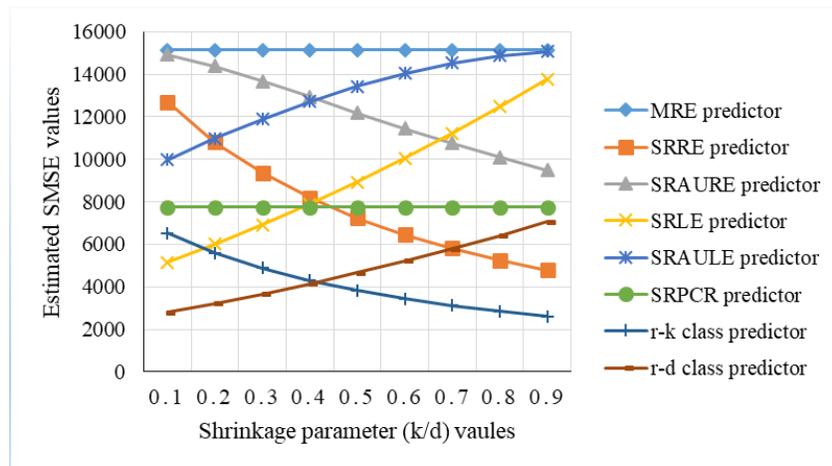

Figure 6. SMSE values of the predictors in the misspecified regression model $((l, p) = (3, 2))$ when n = 50 and $\rho = 0.999$

## 5 Conclusion

Theorem 1 and Theorem 2 give the common form of superiority conditions to compare the estimators (MRE, SRRE, SRAURE, SRLE, SRAULE, SRPCR, SRrk and SRrd) and their respective predictors in MSEM criterion in the misspecified linear regression model when the prior information of the regression coefficients is incomplete, and the multicollinearity exists among the explanatory variables.

From the simulation study, it can be identified the superior estimators and predictors over the others when the conditions are different. The results obtained in this research will produce significant improvements in the parameter estimation in misspecified regression models with incomplete prior information, and the results are applicable to real-world applications.

## Conflicts of Interest

The authors declare that they have no conflicts of interest.

## Appendix A: Lemmas

**Lemma A1:** (Trenkler and Toutenburg, 1990)
Let $\hat{\beta}_1$ and $\hat{\beta}_2$ be two linear estimator of $\beta$. Suppose that $D = D(\hat{\beta}_1) - D(\hat{\beta}_2)$ is positive definite, then $\Delta = MSEM(\hat{\beta}_1) - MSEM(\hat{\beta}_2)$ is non negative if and only if $b_2'(D + b_1 b_1')^{-1} b_2 \leq 1$, where $D(\hat{\beta}_j)$, $MSE(\hat{\beta}_j)$ and $b_j$ denote dispersion matrix, mean square error matrix and bias vector of $\hat{\beta}_j$ respectively, $j = 1, 2$.

**Lemma A2:** (Wang et al., 2006)
Let $n \times n$ matrices $M > 0, N \geq 0$, then $M > N$ if and only if $\lambda_* < 1$, where $\lambda_*$ is the largest eigenvalue of the matrix $NM^{-1}$.

**Lemma A3:** (Baksalary and Kala, 1983)
Let B$\geq 0$ of type $n \times n$ matrix, $b$ is a $n \times 1$ vector and $\lambda$ is a positive real number. Then the following conditions are equivalent.
  i. $\lambda B - bb' \geq 0$
  ii. $B \geq 0$, $b \in \Re(B)$ and $b'B^{-1}b \leq \lambda$, where $\Re(B)$ stands for column space of $B$ and $B^{-1}$ is a independent choice of g-inverse of $B$.

# Appendix B: Tables

Table B1. Bias vector, Dispersion matrix and MSEM of the estimators

| Estimators ($\hat{\gamma}_{(j)}$) | $Bias(\hat{\gamma}_{(j)})$, $D(\hat{\gamma}_{(j)})$ and $MSEM(\hat{\gamma}_{(j)})$ |
|---|---|
| $\hat{\gamma}_{MRE}$ | $Bias(\hat{\gamma}_{MRE}) = \tau A$ <br> $D(\hat{\gamma}_{MRE}) = \sigma^2 \tau$ <br> $MSEM(\hat{\gamma}_{MRE}) = \sigma^2 \tau + (\tau A)(\tau A)'$ |
| $\hat{\gamma}_{SRRE}$ | $Bias(\hat{\gamma}_{SRRE}) = (1+k)^{-1}(\tau A - k\gamma)$ <br> $D(\hat{\gamma}_{SRRE}) = (1+k)^{-2}\sigma^2 \tau$ <br> $MSEM(\hat{\gamma}_{SRRE}) = (1+k)^{-2}(\sigma^2 \tau + (\tau A - k\gamma)(\tau A - k\gamma)')$ |
| $\hat{\gamma}_{SRAURE}$ | $Bias(\hat{\gamma}_{SRAURE}) = (1+k)^{-2}\left((1+2k)\tau A - k^2\gamma\right)$ <br> $D(\hat{\gamma}_{SRAURE}) = (1+k)^{-4}(1+2k)^2 \sigma^2 \tau$ <br> $MSEM(\hat{\gamma}_{SRAURE}) = (1+k)^{-4}\left((1+2k)^2\sigma^2\tau + \left((1+2k)\tau A - k^2\gamma\right)\left((1+2k)\tau A - k^2\gamma\right)'\right)$ |
| $\hat{\gamma}_{SRLE}$ | $Bias(\hat{\gamma}_{SRLE}) = 2^{-1}\left((1+d)\tau A - (1-d)\gamma\right)$ <br> $D(\hat{\gamma}_{SRLE}) = 2^{-2}(1+d)^2 \sigma^2 \tau$ <br> $MSEM(\hat{\gamma}_{SRLE}) = 2^{-2}\left((1+d)^2\sigma^2\tau + \left((1+d)\tau A - (1-d)\gamma\right)\left((1+d)\tau A - (1-d)\gamma\right)'\right)$ |
| $\hat{\gamma}_{SRAULE}$ | $Bias(\hat{\gamma}_{SRAULE}) = 2^{-2}((1+d)(3-d)\tau A - (1-d)^2\gamma)$ <br> $D(\hat{\gamma}_{SRAULE}) = 2^{-4}(1+d)^2(3-d)^2 \sigma^2 \tau$ <br> $MSEM(\hat{\gamma}_{SRAULE}) = 2^{-4}\left((1+d)^2(3-d)^2\sigma^2\tau + \left((1+d)(3-d)\tau A - (1-d)^2\gamma\right)\left((1+d)(3-d)\tau A - (1-d)^2\gamma\right)'\right)$ |
| $\hat{\gamma}_{SRPCR}$ | $Bias(\hat{\gamma}_{SRPCR}) = (T_h T_h' - I)\gamma + T_h T_h' \tau A$ <br> $D(\hat{\gamma}_{SRPCR}) = \sigma^2 T_h T_h' \tau T_h' T_h$ <br> $MSEM(\hat{\gamma}_{SRPCR}) = \sigma^2 T_h T_h' \tau T_h' T_h + \left((T_h T_h' - I)\gamma + T_h T_h' \tau A\right)\left((T_h T_h' - I)\gamma + T_h T_h' \tau A\right)'$ |
| $\hat{\gamma}_{SRrk}$ | $Bias(\hat{\gamma}_{SRrk}) = (1+k)^{-1}\left((T_h T_h' - (1+k)I)\gamma + T_h T_h' \tau A\right)$ <br> $D(\hat{\gamma}_{SRrk}) = (1+k)^{-2}\sigma^2 T_h T_h' \tau T_h' T_h$ <br> $MSEM(\hat{\gamma}_{SRrk}) = (1+k)^{-2}\left(\sigma^2 T_h T_h' \tau T_h' T_h + \left((T_h T_h' - (1+k)I)\gamma + T_h T_h' \tau A\right)\left((T_h T_h' - (1+k)I)\gamma + T_h T_h' \tau A\right)'\right)$ |
| $\hat{\gamma}_{SRrd}$ | $Bias(\hat{\gamma}_{SRrd}) = 2^{-1}(1+d)\left((T_h T_h' - 2(1+d)^{-1}I)\gamma + T_h T_h' \tau A\right)$ <br> $D(\hat{\gamma}_{SRrd}) = 2^{-2}(1+d)^2 \sigma^2 T_h T_h' \tau T_h' T_h$ <br> $MSEM(\hat{\gamma}_{SRrd}) = 2^{-2}(1+d)^2\left(\sigma^2 T_h T_h' \tau T_h' T_h + \left((T_h T_h' - 2(1+d)^{-1}I)\gamma + T_h T_h' \tau A\right)\left((T_h T_h' - 2(1+d)^{-1}I)\gamma + T_h T_h' \tau A\right)'\right)$ |

Table B2. Estimated SMSE values of the estimators

| $k/d$ | 0.1 | 0.2 | 0.3 | 0.4 | 0.5 | 0.6 | 0.7 | 0.8 | 0.9 |
|---|---|---|---|---|---|---|---|---|---|
| $(l, p) = (4, 0)$ | | | | | | | | | |
| $SMSE(\hat{\gamma}_{MRE})$ | **0.119** | **0.119** | **0.119** | **0.119** | **0.119** | **0.119** | **0.119** | **0.119** | **0.119** |
| $SMSE(\hat{\gamma}_{SRRE})$ | 0.545 | 1.672 | 3.174 | 4.860 | 6.621 | 8.388 | 10.124 | 11.805 | 13.420 |
| $SMSE(\hat{\gamma}_{SRAURE})$ | 0.116 | 0.142 | 0.246 | 0.456 | 0.777 | 1.206 | 1.729 | 2.333 | 3.002 |
| $SMSE(\hat{\gamma}_{SRLE})$ | 12.104 | 9.550 | 7.302 | 5.359 | 3.722 | 2.390 | 1.364 | 0.644 | 0.229 |
| $SMSE(\hat{\gamma}_{SRAULE})$ | 2.450 | 1.545 | 0.930 | 0.536 | 0.304 | 0.182 | 0.131 | 0.117 | 0.117 |
| $SMSE(\hat{\gamma}_{SRPCR})$ | 3.336 | 3.336 | 3.336 | 3.336 | 3.336 | 3.336 | 3.336 | 3.336 | 3.336 |
| $SMSE(\hat{\gamma}_{rk})$ | 3.695 | 4.731 | 6.131 | 7.713 | 9.369 | 11.035 | 12.674 | 14.263 | 15.791 |
| $SMSE(\hat{\gamma}_{rd})$ | 14.546 | 12.132 | 10.011 | 8.182 | 6.644 | 5.399 | 4.445 | 3.784 | 3.414 |
| $(l, p) = (3, 1)$ | | | | | | | | | |
| $SMSE(\hat{\gamma}_{MRE})$ | 2.974 | 2.974 | 2.974 | 2.974 | 2.974 | 2.974 | 2.974 | 2.974 | 2.974 |
| $SMSE(\hat{\gamma}_{SRRE})$ | 1.171 | **0.395** | **0.255** | **0.511** | 1.014 | 1.668 | 2.411 | 3.202 | 4.014 |
| $SMSE(\hat{\gamma}_{SRAURE})$ | 2.771 | 2.323 | 1.803 | 1.312 | 0.900 | 0.588 | **0.379** | **0.269** | **0.248** |
| $SMSE(\hat{\gamma}_{SRLE})$ | 3.349 | 2.156 | 1.251 | 0.634 | **0.305** | **0.263** | 0.509 | 1.043 | 1.865 |
| $SMSE(\hat{\gamma}_{SRAULE})$ | **0.259** | 0.437 | 0.767 | 1.185 | 1.633 | 2.064 | 2.440 | 2.730 | 2.912 |
| $SMSE(\hat{\gamma}_{SRPCR})$ | 13.712 | 13.712 | 13.712 | 13.712 | 13.712 | 13.712 | 13.712 | 13.712 | 13.712 |
| $SMSE(\hat{\gamma}_{rk})$ | 13.273 | 13.276 | 13.541 | 13.959 | 14.464 | 15.015 | 15.586 | 16.159 | 16.725 |
| $SMSE(\hat{\gamma}_{rd})$ | 16.263 | 15.395 | 14.673 | 14.097 | 13.667 | 13.384 | 13.247 | 13.256 | 13.411 |
| $(l, p) = (2, 2)$ | | | | | | | | | |
| $SMSE(\hat{\gamma}_{MRE})$ | 8.499 | 8.499 | 8.499 | 8.499 | 8.499 | 8.499 | 8.499 | 8.499 | 8.499 |
| $SMSE(\hat{\gamma}_{SRRE})$ | 4.873 | 2.617 | 1.251 | 0.477 | **0.102** | **0.001** | **0.085** | **0.297** | **0.598** |
| $SMSE(\hat{\gamma}_{SRAURE})$ | 8.128 | 7.285 | 6.253 | 5.197 | 4.203 | 3.314 | 2.544 | 1.897 | 1.365 |
| $SMSE(\hat{\gamma}_{SRLE})$ | **0.347** | **0.040** | **0.036** | **0.336** | 0.938 | 1.844 | 3.053 | 4.565 | 6.381 |
| $SMSE(\hat{\gamma}_{SRAULE})$ | 1.792 | 2.787 | 3.847 | 4.904 | 5.898 | 6.780 | 7.508 | 8.051 | 8.386 |
| $SMSE(\hat{\gamma}_{SRPCR})$ | 12.232 | 12.232 | 12.232 | 12.232 | 12.232 | 12.232 | 12.232 | 12.232 | 12.232 |
| $SMSE(\hat{\gamma}_{rk})$ | 12.916 | 13.542 | 14.111 | 14.627 | 15.095 | 15.521 | 15.910 | 16.265 | 16.591 |
| $SMSE(\hat{\gamma}_{rd})$ | 16.326 | 15.784 | 15.264 | 14.765 | 14.288 | 13.834 | 13.400 | 12.989 | 12.599 |

Table B3. Estimated SMSE values of the predictors

| $k/d$ | 0.1 | 0.2 | 0.3 | 0.4 | 0.5 | 0.6 | 0.7 | 0.8 | 0.9 |
|---|---|---|---|---|---|---|---|---|---|
| $(l, p) = (4, 0)$ | | | | | | | | | |
| $SMSE(\hat{y}_{MRE})$ | 7954.9 | 7954.9 | 7954.9 | 7954.9 | 7954.9 | 7954.9 | 7954.9 | 7954.9 | 7954.9 |
| $SMSE(\hat{y}_{SRRE})$ | 6460.6 | 5334.0 | 4464.9 | 3781.5 | **3235.0** | **2791.9** | **2428.0** | **2125.8** | **1872.5** |
| $SMSE(\hat{y}_{SRAURE})$ | 7812.7 | 7481.8 | 7060.7 | 6606.0 | 6149.6 | 5709.1 | 5293.2 | 4905.9 | 4548.1 |
| $SMSE(\hat{y}_{SRLE})$ | **2076.5** | **2541.7** | **3053.8** | **3613.0** | 4219.3 | 4872.4 | 5572.6 | 6319.7 | 7113.8 |
| $SMSE(\hat{y}_{SRAULE})$ | 4838.6 | 5428.8 | 5977.7 | 6474.8 | 6910.9 | 7278.3 | 7570.6 | 7783.0 | 7911.8 |
| $SMSE(\hat{y}_{SRPCR})$ | 7954.8 | 7954.8 | 7954.8 | 7954.8 | 7954.8 | 7954.8 | 7954.8 | 7954.8 | 7954.8 |
| $SMSE(\hat{y}_{SRrk})$ | 6460.5 | 5333.9 | 4464.8 | 3781.4 | 3235.1 | 2790.0 | 2428.1 | 2125.9 | 1872.6 |
| $SMSE(\hat{y}_{SRrd})$ | 2076.6 | 2541.8 | 3053.9 | 3613.1 | 4219.2 | 4872.3 | 5572.5 | 6319.6 | 7113.7 |
| $(l, p) = (3, 1)$ | | | | | | | | | |
| $SMSE(\hat{y}_{MRE})$ | 3496.8 | 3496.8 | 3496.8 | 3496.8 | 3496.8 | 3496.8 | 3496.8 | 3496.8 | 3496.8 |
| $SMSE(\hat{y}_{SRRE})$ | 2966.2 | 2557.4 | 2235.1 | 1976.2 | **1764.7** | **1589.4** | **1442.4** | **1317.7** | **1210.9** |
| $SMSE(\hat{y}_{SRAURE})$ | 3446.7 | 3330.0 | 3180.7 | 3018.3 | 2854.2 | 2694.5 | 2542.4 | 2399.5 | 2266.3 |
| $SMSE(\hat{y}_{SRLE})$ | **1297.0** | **1488.7** | **1693.5** | **1911.5** | 2142.7 | 2387.1 | 2644.7 | 2915.5 | 3199.5 |
| $SMSE(\hat{y}_{SRAULE})$ | 2374.6 | 2592.1 | 2792.0 | 2971.3 | 3127.3 | 3257.9 | 3361.4 | 3436.3 | 3481.6 |
| $SMSE(\hat{y}_{SRPCR})$ | 5417.8 | 5417.8 | 5417.8 | 5417.8 | 5417.8 | 5417.8 | 5417.8 | 5417.8 | 5417.8 |
| $SMSE(\hat{y}_{SRrk})$ | 4572.5 | 3922.8 | 3412.0 | 3002.4 | 2668.7 | 2392.8 | 2161.8 | 1966.4 | 1799.3 |
| $SMSE(\hat{y}_{SRrd})$ | 1934.1 | 2234.4 | 2556.5 | 2900.2 | 3265.6 | 3652.7 | 4061.5 | 4491.9 | 4944.0 |
| $(l, p) = (2, 2)$ | | | | | | | | | |
| $SMSE(\hat{y}_{MRE})$ | 4864.5 | 4864.5 | 4864.5 | 4864.5 | 4864.5 | 4864.5 | 4864.5 | 4864.5 | 4864.5 |
| $SMSE(\hat{y}_{SRRE})$ | 4110.3 | 3530.3 | 3073.9 | 2707.9 | 2409.5 | 2162.6 | 1955.9 | 1780.8 | 1631.1 |
| $SMSE(\hat{y}_{SRAURE})$ | 4793.3 | 4627.3 | 4415.0 | 4184.3 | 3951.3 | 3724.7 | 3509.1 | 3306.6 | 3118.0 |
| $SMSE(\hat{y}_{SRLE})$ | 1751.9 | 2020.9 | 2309.1 | 2616.5 | 2943.2 | 3289.0 | 3654.1 | 4038.3 | 4441.8 |
| $SMSE(\hat{y}_{SRAULE})$ | 3271.3 | 3579.5 | 3863.0 | 4117.5 | 4339.1 | 4524.8 | 4671.9 | 4778.4 | 4842.9 |
| $SMSE(\hat{y}_{SRPCR})$ | 847.0 | 847.0 | 847.0 | 847.0 | 847.0 | 847.0 | 847.0 | 847.0 | 847.0 |
| $SMSE(\hat{y}_{SRrk})$ | 737.8 | 652.6 | 584.6 | 529.3 | **483.6** | **445.3** | **412.8** | **385.0** | **360.9** |
| $SMSE(\hat{y}_{SRrd})$ | **380.3** | **423.1** | **468.1** | **515.4** | 565.0 | 616.8 | 670.9 | 727.3 | 786.0 |

Table B4. Estimated SMSE values of the estimators when $n = 50$ and $\rho = 0.9$

| k/d | 0.1 | 0.2 | 0.3 | 0.4 | 0.5 | 0.6 | 0.7 | 0.8 | 0.9 |
|---|---|---|---|---|---|---|---|---|---|
| $(l, p) = (3, 0)$ | | | | | | | | | |
| $SMSE(\hat{\gamma}_{MRE})$ | **5.73** | **5.73** | **5.73** | **5.73** | **5.73** | **5.73** | **5.73** | 5.73 | 5.73 |
| $SMSE(\hat{\gamma}_{SRRE})$ | 7.38 | 12.01 | 18.23 | 25.25 | 32.60 | 39.98 | 47.23 | 54.27 | 61.03 |
| $SMSE(\hat{\gamma}_{SRAURE})$ | 5.70 | 5.78 | 6.18 | 7.01 | 8.32 | 10.08 | 12.24 | 14.74 | 17.52 |
| $SMSE(\hat{\gamma}_{SRLE})$ | 55.52 | 44.84 | 35.44 | 27.33 | 20.51 | 14.98 | 10.74 | 7.78 | 6.11 |
| $SMSE(\hat{\gamma}_{SRAULE})$ | 15.23 | 11.48 | 8.95 | 7.34 | 6.40 | 5.93 | 5.74 | **5.70** | **5.71** |
| $SMSE(\hat{\gamma}_{SRPCR})$ | 165.01 | 165.01 | 165.01 | 165.01 | 165.01 | 165.01 | 165.01 | 165.01 | 165.01 |
| $SMSE(\hat{\gamma}_{SRrk})$ | 165.18 | 166.59 | 168.68 | 171.13 | 173.74 | 176.40 | 179.04 | 181.60 | 184.08 |
| $SMSE(\hat{\gamma}_{SRrd})$ | 182.06 | 178.16 | 174.76 | 171.87 | 169.47 | 167.58 | 166.18 | 165.29 | 164.90 |
| $(l, p) = (4, 1)$ | | | | | | | | | |
| $SMSE(\hat{\gamma}_{MRE})$ | 31.497 | 31.497 | 31.497 | 31.497 | 31.497 | 31.497 | 31.497 | 31.497 | 31.497 |
| $SMSE(\hat{\gamma}_{SRRE})$ | 17.883 | 9.773 | 5.207 | 2.969 | **2.280** | **2.631** | **3.679** | **5.192** | 7.009 |
| $SMSE(\hat{\gamma}_{SRAURE})$ | 30.085 | 26.888 | 23.008 | 19.078 | 15.433 | 12.229 | 9.522 | 7.312 | **5.568** |
| $SMSE(\hat{\gamma}_{SRLE})$ | **5.503** | **3.268** | **2.314** | **2.640** | 4.248 | 7.136 | 11.305 | 16.755 | 23.486 |
| $SMSE(\hat{\gamma}_{SRAULE})$ | 6.961 | 10.369 | 14.144 | 17.998 | 21.683 | 24.985 | 27.733 | 29.792 | 31.066 |
| $SMSE(\hat{\gamma}_{SRPCR})$ | 44.136 | 44.136 | 44.136 | 44.136 | 44.136 | 44.136 | 44.136 | 44.136 | 44.136 |
| $SMSE(\hat{\gamma}_{SRrk})$ | 37.801 | 34.469 | 33.032 | 32.809 | 33.368 | 34.431 | 35.814 | 37.393 | 39.085 |
| $SMSE(\hat{\gamma}_{SRrd})$ | 37.694 | 35.326 | 33.729 | 32.903 | 32.848 | 33.563 | 35.050 | 37.308 | 40.336 |
| $(l, p) = (2, 2)$ | | | | | | | | | |
| $SMSE(\hat{\gamma}_{MRE})$ | 61.869 | 61.869 | 61.869 | 61.869 | 61.869 | 61.869 | 61.869 | 61.869 | 61.869 |
| $SMSE(\hat{\gamma}_{SRRE})$ | 41.520 | 27.719 | 18.282 | 11.828 | **7.455** | **4.559** | **2.724** | **1.661** | **1.161** |
| $SMSE(\hat{\gamma}_{SRAURE})$ | 59.848 | 55.213 | 49.447 | 43.407 | 37.559 | 32.139 | 27.250 | 22.918 | 19.130 |
| $SMSE(\hat{\gamma}_{SRLE})$ | **1.533** | **3.238** | **6.192** | **10.397** | 15.851 | 22.555 | 30.509 | 39.712 | 50.166 |
| $SMSE(\hat{\gamma}_{SRAULE})$ | 22.190 | 28.818 | 35.416 | 41.703 | 47.435 | 52.406 | 56.448 | 59.428 | 61.254 |
| $SMSE(\hat{\gamma}_{SRPCR})$ | 38.090 | 38.090 | 38.090 | 38.090 | 38.090 | 38.090 | 38.090 | 38.090 | 38.090 |
| $SMSE(\hat{\gamma}_{SRrk})$ | 34.346 | 32.175 | 31.013 | 30.508 | 30.438 | 30.657 | 31.066 | 31.601 | 32.216 |
| $SMSE(\hat{\gamma}_{SRrd})$ | 31.708 | 30.913 | 30.494 | 30.451 | 30.784 | 31.493 | 32.579 | 34.040 | 35.877 |

Table B5. Estimated SMSE values of the estimators when $n = 50$ and $\rho = 0.99$

| k/d | 0.1 | 0.2 | 0.3 | 0.4 | 0.5 | 0.6 | 0.7 | 0.8 | 0.9 |
|---|---|---|---|---|---|---|---|---|---|
| $(l, p) = (5, 0)$ | | | | | | | | | |
| $SMSE(\hat{\gamma}_{MRE})$ | 5.89 | 5.89 | **5.89** | **5.89** | **5.89** | **5.89** | 5.89 | 5.89 | 5.89 |
| $SMSE(\hat{\gamma}_{SRRE})$ | 7.38 | 11.93 | 18.13 | 25.16 | 32.54 | 39.97 | 47.29 | 54.38 | 61.21 |
| $SMSE(\hat{\gamma}_{SRAURE})$ | **5.84** | **5.88** | 6.23 | 7.02 | 8.29 | 10.02 | 12.16 | 14.65 | 17.42 |
| $SMSE(\hat{\gamma}_{SRLE})$ | 55.65 | 44.87 | 35.41 | 27.25 | 20.41 | 14.88 | 10.67 | 7.76 | 6.17 |
| $SMSE(\hat{\gamma}_{SRAULE})$ | 15.13 | 11.41 | 8.91 | 7.34 | 6.44 | 6.01 | **5.85** | **5.84** | **5.87** |
| $SMSE(\hat{\gamma}_{SRPCR})$ | 195.79 | 195.79 | 195.79 | 195.79 | 195.79 | 195.79 | 195.79 | 195.79 | 195.79 |
| $SMSE(\hat{\gamma}_{SRrk})$ | 195.66 | 196.47 | 197.82 | 199.45 | 201.22 | 203.04 | 204.85 | 206.63 | 208.36 |
| $SMSE(\hat{\gamma}_{SRrd})$ | 206.95 | 204.25 | 201.91 | 199.94 | 198.34 | 197.10 | 196.22 | 195.71 | 195.57 |
| $(l, p) = (4, 1)$ | | | | | | | | | |
| $SMSE(\hat{\gamma}_{MRE})$ | 19.939 | 19.939 | 19.939 | 19.939 | 19.939 | 19.939 | 19.939 | 19.939 | 19.939 |
| $SMSE(\hat{\gamma}_{SRRE})$ | 10.064 | **5.145** | **3.333** | **3.494** | 4.914 | 7.131 | 9.844 | 12.850 | 16.014 |
| $SMSE(\hat{\gamma}_{SRAURE})$ | 18.862 | 16.462 | 13.629 | 10.873 | 8.458 | 6.497 | **5.017** | **4.004** | **3.417** |
| $SMSE(\hat{\gamma}_{SRLE})$ | 13.416 | 8.899 | 5.691 | 3.795 | **3.209** | **3.934** | 5.969 | 9.315 | 13.972 |
| $SMSE(\hat{\gamma}_{SRAULE})$ | **3.866** | 5.457 | 7.647 | 10.142 | 12.685 | 15.060 | 17.091 | 18.640 | 19.609 |
| $SMSE(\hat{\gamma}_{SRPCR})$ | 47.194 | 47.194 | 47.194 | 47.194 | 47.194 | 47.194 | 47.194 | 47.194 | 47.194 |
| $SMSE(\hat{\gamma}_{SRrk})$ | 42.791 | 41.218 | 41.373 | 42.592 | 44.459 | 46.710 | 49.175 | 51.742 | 54.339 |
| $SMSE(\hat{\gamma}_{SRrd})$ | 52.213 | 48.337 | 45.290 | 43.073 | 41.686 | 41.128 | 41.400 | 42.502 | 44.433 |
| $(l, p) = (3, 2)$ | | | | | | | | | |
| $SMSE(\hat{\gamma}_{MRE})$ | 49.955 | 49.955 | 49.955 | 49.955 | 49.955 | 49.955 | 49.955 | 49.955 | 49.955 |
| $SMSE(\hat{\gamma}_{SRRE})$ | 31.645 | 19.652 | 11.823 | 6.804 | **3.718** | **1.979** | **1.190** | **1.076** | **1.443** |
| $SMSE(\hat{\gamma}_{SRAURE})$ | 48.113 | 43.906 | 38.710 | 33.317 | 28.157 | 23.440 | 19.254 | 15.617 | 12.509 |
| $SMSE(\hat{\gamma}_{SRLE})$ | **1.112** | **1.367** | **2.914** | **5.755** | 9.889 | 15.316 | 22.036 | 30.049 | 39.355 |
| $SMSE(\hat{\gamma}_{SRAULE})$ | 15.014 | 20.589 | 26.283 | 31.807 | 36.908 | 41.372 | 45.025 | 47.731 | 49.394 |
| $SMSE(\hat{\gamma}_{SRPCR})$ | 44.297 | 44.297 | 44.297 | 44.297 | 44.297 | 44.297 | 44.297 | 44.297 | 44.297 |
| $SMSE(\hat{\gamma}_{SRrk})$ | 40.730 | 38.850 | 38.034 | 37.900 | 38.206 | 38.796 | 39.565 | 40.445 | 41.389 |
| $SMSE(\hat{\gamma}_{SRrd})$ | 40.613 | 39.293 | 38.406 | 37.951 | 37.928 | 38.337 | 39.179 | 40.452 | 42.158 |

Table B6. Estimated SMSE values of the estimators when $n = 50$ and $\rho = 0.999$

| k/d | 0.1 | 0.2 | 0.3 | 0.4 | 0.5 | 0.6 | 0.7 | 0.8 | 0.9 |
|---|---|---|---|---|---|---|---|---|---|
| $(l, p) = (5, 0)$ | | | | | | | | | |
| $SMSE(\hat{\gamma}_{MRE})$ | **2.39** | **2.39** | **2.39** | **2.39** | **2.39** | **2.39** | **2.39** | **2.39** | **2.39** |
| $SMSE(\hat{\gamma}_{SRRE})$ | 4.52 | 9.46 | 15.90 | 23.06 | 30.49 | 37.93 | 45.21 | 52.26 | 59.02 |
| $SMSE(\hat{\gamma}_{SRAURE})$ | 2.41 | 2.60 | 3.13 | 4.11 | 5.55 | 7.44 | 9.71 | 12.30 | 15.16 |
| $SMSE(\hat{\gamma}_{SRLE})$ | 53.51 | 42.81 | 33.36 | 25.17 | 18.23 | 12.55 | 8.13 | 4.96 | 3.05 |
| $SMSE(\hat{\gamma}_{SRAULE})$ | 12.80 | 8.91 | 6.23 | 4.47 | 3.41 | 2.81 | 2.53 | 2.42 | 2.40 |
| $SMSE(\hat{\gamma}_{SRPCR})$ | 196.31 | 196.31 | 196.31 | 196.31 | 196.31 | 196.31 | 196.31 | 196.31 | 196.31 |
| $SMSE(\hat{\gamma}_{SRrk})$ | 196.66 | 197.70 | 199.12 | 200.72 | 202.39 | 204.08 | 205.74 | 207.35 | 208.90 |
| $SMSE(\hat{\gamma}_{SRrd})$ | 207.64 | 205.19 | 203.05 | 201.19 | 199.64 | 198.38 | 197.41 | 196.75 | 196.38 |
| $(l, p) = (4, 1)$ | | | | | | | | | |
| $SMSE(\hat{\gamma}_{MRE})$ | 13.160 | 13.160 | 13.160 | 13.160 | 13.160 | 13.160 | 13.160 | 13.160 | 13.160 |
| $SMSE(\hat{\gamma}_{SRRE})$ | 5.343 | **2.031** | **1.501** | **2.706** | 4.988 | 7.929 | 11.255 | 14.786 | 18.405 |
| $SMSE(\hat{\gamma}_{SRAURE})$ | 12.276 | 10.327 | 8.073 | 5.949 | 4.176 | 2.841 | **1.962** | **1.516** | **1.461** |
| $SMSE(\hat{\gamma}_{SRLE})$ | 15.440 | 10.116 | 6.059 | 3.270 | **1.749** | **1.496** | 2.510 | 4.792 | 8.342 |
| $SMSE(\hat{\gamma}_{SRAULE})$ | **1.478** | 2.206 | 3.608 | 5.401 | 7.336 | 9.204 | 10.835 | 12.095 | 12.889 |
| $SMSE(\hat{\gamma}_{SRPCR})$ | 45.790 | 45.790 | 45.790 | 45.790 | 45.790 | 45.790 | 45.790 | 45.790 | 45.790 |
| $SMSE(\hat{\gamma}_{SRrk})$ | 42.368 | 41.593 | 42.412 | 44.190 | 46.535 | 49.198 | 52.022 | 54.906 | 57.784 |
| $SMSE(\hat{\gamma}_{SRrd})$ | 55.431 | 51.070 | 47.532 | 44.815 | 42.922 | 41.850 | 41.602 | 42.175 | 43.571 |
| $(l, p) = (3, 2)$ | | | | | | | | | |
| $SMSE(\hat{\gamma}_{MRE})$ | 42.748 | 42.748 | 42.748 | 42.748 | 42.748 | 42.748 | 42.748 | 42.748 | 42.748 |
| $SMSE(\hat{\gamma}_{SRRE})$ | 26.110 | 15.416 | 8.618 | 4.435 | **2.035** | **0.870** | **0.566** | **0.866** | **1.589** |
| $SMSE(\hat{\gamma}_{SRAURE})$ | 41.064 | 37.224 | 32.498 | 27.618 | 22.976 | 18.767 | 15.066 | 11.885 | 9.204 |
| $SMSE(\hat{\gamma}_{SRLE})$ | **0.971** | **0.590** | **1.465** | **3.595** | 6.982 | 11.624 | 17.521 | 24.675 | 33.084 |
| $SMSE(\hat{\gamma}_{SRAULE})$ | 11.362 | 16.241 | 21.300 | 26.256 | 30.864 | 34.916 | 38.244 | 40.715 | 42.235 |
| $SMSE(\hat{\gamma}_{SRPCR})$ | 44.641 | 44.641 | 44.641 | 44.641 | 44.641 | 44.641 | 44.641 | 44.641 | 44.641 |
| $SMSE(\hat{\gamma}_{SRrk})$ | 41.468 | 39.910 | 39.365 | 39.461 | 39.965 | 40.726 | 41.646 | 42.659 | 43.722 |
| $SMSE(\hat{\gamma}_{SRrd})$ | 42.850 | 41.326 | 40.233 | 39.571 | 39.339 | 39.538 | 40.168 | 41.229 | 42.720 |

Table B7. Estimated SMSE values of the predictors when $n = 50$ and $\rho = 0.9$

| k/d | 0.1 | 0.2 | 0.3 | 0.4 | 0.5 | 0.6 | 0.7 | 0.8 | 0.9 |
|---|---|---|---|---|---|---|---|---|---|
| $(l, p) = (5, 0)$ | | | | | | | | | |
| $SMSE(\hat{y}_{MRE})$ | 26713 | 26713 | 26713 | 26713 | 26713 | 26713 | 26713 | 26713 | 26713 |
| $SMSE(\hat{y}_{SRRE})$ | 22354 | 19021 | 16413 | 14333 | 12646 | 11258 | 10101 | 9127 | 8298 |
| $SMSE(\hat{y}_{SRAURE})$ | 26300 | 25340 | 24112 | 22781 | 21439 | 20136 | 18900 | 17741 | 16665 |
| $SMSE(\hat{y}_{SRLE})$ | 8966 | 10464 | 12081 | 13816 | 15669 | 17641 | 19731 | 21940 | 24267 |
| $SMSE(\hat{y}_{SRAULE})$ | 17540 | 19304 | 20931 | 22396 | 23674 | 24747 | 25598 | 26214 | 26588 |
| $SMSE(\hat{y}_{SRPCR})$ | 26257 | 26257 | 26257 | 26257 | 26257 | 26257 | 26257 | 26257 | 26257 |
| $SMSE(\hat{y}_{SRrk})$ | 21978 | 18705 | 16144 | 14101 | **12444** | **11080** | **9944** | **8987** | **8172** |
| $SMSE(\hat{y}_{SRrd})$ | **8829** | **10301** | **11889** | **13593** | 15413 | 17350 | 19402 | 21571 | 23856 |
| $(l, p) = (4, 1)$ | | | | | | | | | |
| $SMSE(\hat{y}_{MRE})$ | 16251 | 16251 | 16251 | 16251 | 16251 | 16251 | 16251 | 16251 | 16251 |
| $SMSE(\hat{y}_{SRRE})$ | 13669 | 11690 | 10139 | 8898 | 7891 | 7060 | 6366 | 5781 | 5282 |
| $SMSE(\hat{y}_{SRAURE})$ | 16007 | 15438 | 14711 | 13922 | 13126 | 12352 | 11618 | 10929 | 10288 |
| $SMSE(\hat{y}_{SRLE})$ | 5684 | 6584 | 7553 | 8590 | 9695 | 10869 | 12112 | 13423 | 14803 |
| $SMSE(\hat{y}_{SRAULE})$ | 10809 | 11858 | 12824 | 13693 | 14451 | 15087 | 15591 | 15956 | 16177 |
| $SMSE(\hat{y}_{SRPCR})$ | 8636 | 8636 | 8636 | 8636 | 8636 | 8636 | 8636 | 8636 | 8636 |
| $SMSE(\hat{y}_{SRrk})$ | 7314 | 6299 | 5502 | 4863 | **4342** | **3912** | **3552** | **3248** | **2988** |
| $SMSE(\hat{y}_{SRrd})$ | **3198** | **3666** | **4167** | **4703** | 5273 | 5877 | 6516 | 7188 | 7895 |
| $(l, p) = (3, 2)$ | | | | | | | | | |
| $SMSE(\hat{y}_{MRE})$ | 14936 | 14936 | 14936 | 14936 | 14936 | 14936 | 14936 | 14936 | 14936 |
| $SMSE(\hat{y}_{SRRE})$ | 12511 | 10657 | 9207 | 8051 | 7114 | 6344 | 5702 | 5161 | 4702 |
| $SMSE(\hat{y}_{SRAURE})$ | 14707 | 14172 | 13489 | 12748 | 12002 | 11277 | 10590 | 9945 | 9347 |
| $SMSE(\hat{y}_{SRLE})$ | 5072 | 5903 | 6800 | 7764 | 8794 | 9890 | 11052 | 12280 | 13575 |
| $SMSE(\hat{y}_{SRAULE})$ | 9833 | 10814 | 11719 | 12534 | 13245 | 13842 | 14316 | 14659 | 14867 |
| $SMSE(\hat{y}_{SRPCR})$ | 5972 | 5972 | 5972 | 5972 | 5972 | 5972 | 5972 | 5972 | 5972 |
| $SMSE(\hat{y}_{SRrk})$ | 5054 | 4350 | 3799 | 3358 | **3000** | **2705** | **2458** | **2250** | **2073** |
| $SMSE(\hat{y}_{SRrd})$ | **2216** | **2536** | **2880** | **3248** | 3641 | 4059 | 4500 | 4966 | 5457 |

Table B8. Estimated SMSE values of the predictors when $n = 50$ and $\rho = 0.99$

| k/d | 0.1 | 0.2 | 0.3 | 0.4 | 0.5 | 0.6 | 0.7 | 0.8 | 0.9 |
|---|---|---|---|---|---|---|---|---|---|
| $(l, p) = (5, 0)$ | | | | | | | | | |
| $SMSE(\hat{y}_{MRE})$ | 20010 | 20010 | 20010 | 20010 | 20010 | 20010 | 20010 | 20010 | 20010 |
| $SMSE(\hat{y}_{SRRE})$ | 16798 | 14337 | 12410 | 10871 | 9620 | 8591 | 7731 | 7007 | 6390 |
| $SMSE(\hat{y}_{SRAURE})$ | 19707 | 18999 | 18094 | 17112 | 16122 | 15161 | 14248 | 13392 | 12596 |
| $SMSE(\hat{y}_{SRLE})$ | 6887 | 8001 | 9201 | 10487 | 11859 | 13318 | 14862 | 16492 | 18208 |
| $SMSE(\hat{y}_{SRAULE})$ | 13243 | 14546 | 15748 | 16828 | 17771 | 18562 | 19189 | 19643 | 19918 |
| $SMSE(\hat{y}_{SRPCR})$ | 19969 | 19969 | 19969 | 19969 | 19969 | 19969 | 19969 | 19969 | 19969 |
| $SMSE(\hat{y}_{SRrk})$ | 16763 | 14309 | 12386 | 10849 | **9602** | **8574** | **7717** | **6994** | **6378** |
| $SMSE(\hat{y}_{SRrd})$ | **6875** | **7986** | **9184** | **10467** | 11836 | 13291 | 14832 | 16458 | 18171 |
| $(l, p) = (4, 1)$ | | | | | | | | | |
| $SMSE(\hat{y}_{MRE})$ | 17211 | 17211 | 17211 | 17211 | 17211 | 17211 | 17211 | 17211 | 17211 |
| $SMSE(\hat{y}_{SRRE})$ | 14481 | 12389 | 10748 | 9436 | 8369 | 7490 | 6756 | 6135 | 5607 |
| $SMSE(\hat{y}_{SRAURE})$ | 16953 | 16351 | 15583 | 14749 | 13907 | 13089 | 12313 | 11584 | 10906 |
| $SMSE(\hat{y}_{SRLE})$ | 6033 | 6986 | 8012 | 9109 | 10279 | 11521 | 12835 | 14221 | 15680 |
| $SMSE(\hat{y}_{SRAULE})$ | 11457 | 12566 | 13589 | 14507 | 15308 | 15980 | 16513 | 16899 | 17132 |
| $SMSE(\hat{y}_{SRPCR})$ | 10589 | 10589 | 10589 | 10589 | 10589 | 10589 | 10589 | 10589 | 10589 |
| $SMSE(\hat{y}_{SRrk})$ | 8960 | 7709 | 6725 | 5937 | **5295** | **4765** | **4321** | **3946** | **3625** |
| $SMSE(\hat{y}_{SRrd})$ | **3884** | **4461** | **5080** | **5741** | 6444 | 7189 | 7976 | 8805 | 9676 |
| $(l, p) = (3, 2)$ | | | | | | | | | |
| $SMSE(\hat{y}_{MRE})$ | 15859 | 15859 | 15859 | 15859 | 15859 | 15859 | 15859 | 15859 | 15859 |
| $SMSE(\hat{y}_{SRRE})$ | 13285 | 11318 | 9779 | 8552 | 7556 | 6738 | 6056 | 5482 | 4994 |
| $SMSE(\hat{y}_{SRAURE})$ | 15615 | 15048 | 14323 | 13537 | 12745 | 11976 | 11246 | 10562 | 9927 |
| $SMSE(\hat{y}_{SRLE})$ | 5388 | 6270 | 7223 | 8246 | 9340 | 10503 | 11737 | 13040 | 14414 |
| $SMSE(\hat{y}_{SRAULE})$ | 10443 | 11484 | 12445 | 13310 | 14064 | 14698 | 15200 | 15564 | 15785 |
| $SMSE(\hat{y}_{SRPCR})$ | 7664 | 7664 | 7664 | 7664 | 7664 | 7664 | 7664 | 7664 | 7664 |
| $SMSE(\hat{y}_{SRrk})$ | 6471 | 5557 | 4840 | 4267 | **3802** | **3418** | **3098** | **2828** | **2597** |
| $SMSE(\hat{y}_{SRrd})$ | **2783** | **3198** | **3646** | **4124** | 4635 | 5177 | 5751 | 6357 | 6995 |

Table B9. Estimated SMSE values of the predictors when $n = 50$ and $\rho = 0.999$

| k/d | 0.1 | 0.2 | 0.3 | 0.4 | 0.5 | 0.6 | 0.7 | 0.8 | 0.9 |
|---|---|---|---|---|---|---|---|---|---|
| $(l, p) = (5, 0)$ | | | | | | | | | |
| $SMSE(\hat{y}_{MRE})$ | 16692 | 16692 | 16692 | 16692 | 16692 | 16692 | 16692 | 16692 | 16692 |
| $SMSE(\hat{y}_{SRRE})$ | 14034 | 11998 | 10401 | 9125 | 8088 | 7234 | 6520 | 5918 | 5405 |
| $SMSE(\hat{y}_{SRAURE})$ | 16441 | 15855 | 15107 | 14295 | 13475 | 12679 | 11923 | 11214 | 10555 |
| $SMSE(\hat{y}_{SRLE})$ | 5819 | 6744 | 7741 | 8807 | 9945 | 11153 | 12432 | 13781 | 15201 |
| $SMSE(\hat{y}_{SRAULE})$ | 11091 | 12170 | 13165 | 14060 | 14840 | 15494 | 16013 | 16388 | 16616 |
| $SMSE(\hat{y}_{SRPCR})$ | 16689 | 16689 | 16689 | 16689 | 16689 | 16689 | 16689 | 16689 | 16689 |
| $SMSE(\hat{y}_{SRrk})$ | 14031 | 11995 | 10399 | 9123 | **8087** | **7232** | **6519** | **5917** | **5404** |
| $SMSE(\hat{y}_{SRrd})$ | **5818** | **6743** | **7739** | **8805** | 9943 | 11151 | 12429 | 13778 | 15198 |
| $(l, p) = (4, 1)$ | | | | | | | | | |
| $SMSE(\hat{y}_{MRE})$ | 16427 | 16427 | 16427 | 16427 | 16427 | 16427 | 16427 | 16427 | 16427 |
| $SMSE(\hat{y}_{SRRE})$ | 13834 | 11845 | 10284 | 9036 | 8021 | 7184 | 6484 | 5893 | 5389 |
| $SMSE(\hat{y}_{SRAURE})$ | 16182 | 15610 | 14880 | 14088 | 13288 | 12511 | 11772 | 11080 | 10435 |
| $SMSE(\hat{y}_{SRLE})$ | 5795 | 6704 | 7680 | 8725 | 9838 | 11019 | 12269 | 13587 | 14973 |
| $SMSE(\hat{y}_{SRAULE})$ | 10959 | 12014 | 12985 | 13858 | 14620 | 15258 | 15764 | 16130 | 16352 |
| $SMSE(\hat{y}_{SRPCR})$ | 10765 | 10765 | 10765 | 10765 | 10765 | 10765 | 10765 | 10765 | 10765 |
| $SMSE(\hat{y}_{SRrk})$ | 9110 | 7839 | 6840 | 6039 | **5387** | **4848** | **4396** | **4015** | **3689** |
| $SMSE(\hat{y}_{SRrd})$ | **3952** | **4538** | **5168** | **5840** | 6554 | 7311 | 8111 | 8953 | 9837 |
| $(l, p) = (3, 2)$ | | | | | | | | | |
| $SMSE(\hat{y}_{MRE})$ | 15149 | 15149 | 15149 | 15149 | 15149 | 15149 | 15149 | 15149 | 15149 |
| $SMSE(\hat{y}_{SRRE})$ | 12696 | 10820 | 9352 | 8182 | 7233 | 6452 | 5801 | 5253 | 4787 |
| $SMSE(\hat{y}_{SRAURE})$ | 14917 | 14376 | 13685 | 12936 | 12181 | 11447 | 10752 | 10100 | 9494 |
| $SMSE(\hat{y}_{SRLE})$ | 5163 | 6006 | 6915 | 7891 | 8934 | 10043 | 11219 | 12463 | 13772 |
| $SMSE(\hat{y}_{SRAULE})$ | 9986 | 10979 | 11895 | 12719 | 13439 | 14042 | 14521 | 14868 | 15079 |
| $SMSE(\hat{y}_{SRPCR})$ | 7736 | 7736 | 7736 | 7736 | 7736 | 7736 | 7736 | 7736 | 7736 |
| $SMSE(\hat{y}_{SRrk})$ | 6527 | 5600 | 4874 | 4294 | **3823** | **3435** | **3111** | **2838** | **2605** |
| $SMSE(\hat{y}_{SRrd})$ | **2793** | **3213** | **3665** | **4150** | 4667 | 5216 | 5798 | 6412 | 7058 |

# Appendix C: Corollaries

**Corollary C1:**

a) $\hat{\gamma}_{SRRE}$ is superior to $\hat{\gamma}_{MRE}$ in MSEM sense when the regression model is misspecified due to excluding relevant variables if and only if

$$(\tau A - k\gamma)'(k(2+k)\sigma^2\tau + (1+k)^2(\tau A)(\tau A)')^{-1}(\tau A - k\gamma) \leq 1$$

**Proof:** Consider $D_{(i,j)} = D(\hat{\gamma}_{MRE}) - D(\hat{\gamma}_{SRRE}) = \sigma^2\tau - (1+k)^{-2}\sigma^2\tau$
$$= ((1+k)^2 - 1)(1+k)^{-2}\sigma^2\tau$$
$$= k(2+k)(1+k)^{-2}\sigma^2\tau$$

Since $k > 0$ and $\tau > 0$, hence $D_{(i,j)} > 0$. This completes the proof.

b) If $A \geq 0$, $\hat{y}_{SRRE}$ is superior to $\hat{y}_{MRE}$ in MSEM sense when the regression model is misspecified due to excluding relevant variables if and only if $\theta \in \mathcal{R}(A)$ and $\theta'A^{-1}\theta \leq 1$, where $A = X_*(k(2+k)(1+k)^{-2}\sigma^2\tau + (\tau A)(\tau A)' - (1+k)^{-2}(\tau A - k\gamma)(\tau A - k\gamma)' + k^2(1+k)^{-2}(\gamma + \tau A)(\gamma + \tau A)')X'_* + \delta\delta'$, $\mathcal{R}(A)$ stands for column space of $A$ and $A^{-1}$ is an independent choice of g-inverse of $A$ and $\theta = \delta + k(1+k)^{-1}X_*(\gamma + \tau A)$.

**Corollary C2:**

a) $\hat{\gamma}_{SRAURE}$ is superior to $\hat{\gamma}_{MRE}$ in MSEM sense when the regression model is misspecified due to excluding relevant variables if and only if

$$\big((1+2k)\tau A - k^2\gamma\big)'(k^2(k^2+4k+2)\sigma^2\tau + (1+k)^4(\tau A)(\tau A)')^{-1}\big((1+2k)\tau A - k^2\gamma\big) \leq 1$$

**Proof:** Consider $D_{(i,j)} = D(\hat{\gamma}_{MRE}) - D(\hat{\gamma}_{SRAURE}) = \sigma^2\tau - (1+k)^{-4}(1+2k)^2\sigma^2\tau$
$$= ((1+k)^4 - (1+2k)^2)(1+k)^{-4}\sigma^2\tau$$
$$= k^2(k^2 + 4k + 2)(1+k)^{-4}\sigma^2\tau$$

Since $k > 0$ and $\tau > 0$, hence $D_{(i,j)} > 0$. This completes the proof.

b) If $A \geq 0$, $\hat{y}_{SRAURE}$ is superior to $\hat{y}_{MRE}$ in MSEM sense when the regression model is misspecified due to excluding relevant variables if and only if $\theta \in \mathcal{R}(A)$ and $\theta'A^{-1}\theta \leq 1$, where $A = X_*\big(k^2(k^2 + 4k + 2)(1+k)^{-4}\sigma^2\tau + (\tau A)(\tau A)' - (1+k)^{-4}\big((1+2k)\tau A - k^2\gamma\big)\big((1+2k)\tau A - k^2\gamma\big)' + k^4(1+k)^{-4}(\gamma + \tau A)(\gamma + \tau A)'\big)X'_* + \delta\delta'$, $\mathcal{R}(A)$ stands for column space of $A$ and $A^{-1}$ is an independent choice of g-inverse of $A$ and $\theta = \delta + k^2(1+k)^{-2}X_*(\gamma + \tau A)$.

**Corollary C3:**

a) $\hat{\gamma}_{SRLE}$ is superior to $\hat{\gamma}_{MRE}$ in MSEM sense when the regression model is misspecified due to excluding relevant variables if and only if

$$\big((1+d)\tau A - (1-d)\gamma\big)'\big((3+d)(1-d)\sigma^2\tau + 2^2(\tau A)(\tau A)'\big)^{-1}\big((1+d)\tau A - (1-d)\gamma\big) \leq 1$$

**Proof:** Consider $D_{(i,j)} = D(\hat{\gamma}_{MRE}) - D(\hat{\gamma}_{SRLE}) = \sigma^2\tau - 2^{-2}(1+d)^2\sigma^2\tau$
$$= (4 - (1+d)^2)2^{-2}\sigma^2\tau$$
$$= 2^{-2}(3+d)(1-d)\sigma^2\tau$$

Since $0 < d < 1$ and $\tau > 0$, hence $D_{(i,j)} > 0$. This completes the proof.

b) If $A \geq 0$, $\hat{y}_{SRLE}$ is superior to $\hat{y}_{MRE}$ in MSEM sense when the regression model is misspecified due to excluding relevant variables if and only if $\theta \in \mathfrak{R}(A)$ and $\theta' A^{-1}\theta \leq 1$, where $A = X_*\big(2^{-2}(3 + d)(1 - d)\sigma^2\tau + (\tau A)(\tau A)' - 2^{-2}\big((1 + d)\tau A - (1 - d)\gamma\big)\big((1 + d)\tau A - (1 - d)\gamma\big)' + 2^{-2}(1 - d)^2(\gamma + \tau A)(\gamma + \tau A)'\big)X_*' + \delta\delta'$, $\mathfrak{R}(A)$ stands for column space of $A$ and $A^{-1}$ is an independent choice of g-inverse of $A$ and $\theta = \delta + 2^{-1}(1 - d) X_*(\gamma + \tau A)$.

**Corollary C4:**

a) $\hat{\gamma}_{SRAULE}$ is superior to $\hat{\gamma}_{MRE}$ in MSEM sense when the regression model is misspecified due to excluding relevant variables if and only if

$$\big((1 + d)(3 - d)\tau A - (1 - d)^2\gamma\big)'\big((7 + 2d - d^2)(1 - d)^2\sigma^2\tau + 2^4(\tau A)(\tau A)'\big)^{-1}\big((1 + d)(3 - d)\tau A - (1 - d)^2\gamma\big) \leq 1$$

**Proof:** Consider $D_{(i,j)} = D(\hat{\gamma}_{MRE}) - D(\hat{\gamma}_{SRAULE}) = \sigma^2\tau - 2^{-4}(1 + d)^2(3 - d)^2\sigma^2\tau$
$$= (2^{-4} - (1 + d)^2(3 - d)^2)2^{-4}\sigma^2\tau$$
$$= 2^{-4}(7 + 2d - d^2)(1 - d)^2\sigma^2\tau$$

Since $0 < d < 1$ and $\tau > 0$, hence $D_{(i,j)} > 0$. This completes the proof.

b) If $A \geq 0$, $\hat{y}_{SRAULE}$ is superior to $\hat{y}_{MRE}$ in MSEM sense when the regression model is misspecified due to excluding relevant variables if and only if $\theta \in \mathfrak{R}(A)$ and $\theta' A^{-1}\theta \leq 1$, where $A = X_*(2^{-4}(7 + 2d - d^2)(1 - d)^2\sigma^2\tau + (\tau A)(\tau A)' - 2^{-4}((1 + d)(3 - d)\tau A - (1 - d)^2\gamma)((1 + d)(3 - d)\tau A - (1 - d)^2\gamma)' + 2^{-4}(1 - d)^4(\gamma + \tau A)(\gamma + \tau A)')X_*' + \delta\delta'$, $\mathfrak{R}(A)$ stands for column space of $A$ and $A^{-1}$ is an independent choice of g-inverse of $A$ and $\theta = \delta + 2^{-2}(1 - d)^2 X_*(\gamma + \tau A)$.

**Corollary C5:**

a) If $\lambda_* < 1$, $\hat{\gamma}_{SRPCR}$ is superior to $\hat{\gamma}_{MRE}$ in MSEM sense when the regression model is misspecified due to excluding relevant variables if and only if

$$\big((T_h T_h' - I)\gamma + T_h T_h' \tau A\big)'\big(\sigma^2(\tau - T_h T_h' \tau T_h' T_h) + (\tau A)(\tau A)'\big)^{-1}\big((T_h T_h' - I)\gamma + T_h T_h' \tau A\big) \leq 1$$

where $\lambda_*$ is the largest eigenvalue of $T_h T_h' \tau T_h' T_h \tau^{-1}$.

**Proof:** Consider $D_{(i,j)} = D(\hat{\gamma}_{MRE}) - D(\hat{\gamma}_{SRPCR}) = \sigma^2\tau - \sigma^2 T_h T_h' \tau T_h' T_h$
$$= \sigma^2(\tau - T_h T_h' \tau T_h' T_h)$$

Since $\tau > 0$, according to Lemma A2 (see Appendix A) $D_{(i,j)} > 0$ if $\lambda_* < 1$, where $\lambda_*$ is the largest eigenvalue of $T_h T_h' \tau T_h' T_h \tau^{-1}$. This completes the proof.

b) If $A \geq 0$, $\hat{y}_{SRPCR}$ is superior to $\hat{y}_{MRE}$ in MSEM sense when the regression model is misspecified due to excluding relevant variables if and only if $\theta \in \mathfrak{R}(A)$ and $\theta' A^{-1}\theta \leq 1$, where $A = X_*\big(\sigma^2(\tau - T_h T_h' \tau T_h' T_h) + (\tau A)(\tau A)' - \big((T_h T_h' - I)\gamma + T_h T_h' \tau A\big)\big((T_h T_h' - I)\gamma + T_h T_h' \tau A\big)' + (I - T_h T_h')(\gamma + \tau A)(\gamma + \tau A)'(I - T_h T_h')'\big)X_*' + \delta\delta'$, $\mathfrak{R}(A)$ stands for column space of $A$ and $A^{-1}$ is an independent choice of g-inverse of $A$ and $\theta = \delta + X_*(I - T_h T_h')(\gamma + \tau A)$.

**Corollary C6:**

a) If $\lambda_* < 1$, $\hat{\gamma}_{SRrk}$ is superior to $\hat{\gamma}_{MRE}$ in MSEM sense when the regression model is misspecified due to excluding relevant variables if and only if $(1+k)^{-2}\big((T_hT_h' - (1+k)I)\gamma + T_hT_h'\tau A\big)'(\sigma^2(\tau - (1+k)^{-2}T_hT_h'\tau T_h'T_h) + (\tau A)(\tau A)')^{-1}\big((T_hT_h' - (1+k)I)\gamma + T_hT_h'\tau A\big) \leq 1$ where $\lambda_*$ is the largest eigenvalue of $(1+k)^{-2}T_hT_h'\tau T_h'T_h\tau^{-1}$.

**Proof:** Consider $D_{(i,j)} = D(\hat{\gamma}_{MRE}) - D(\hat{\gamma}_{SRrk}) = \sigma^2\tau - (1+k)^{-2}\sigma^2 T_hT_h'\tau T_h'T_h$
$= \sigma^2(\tau - (1+k)^{-2}T_hT_h'\tau T_h'T_h)$

Since $\tau > 0$, according to Lemma A2 (see Appendix A) $D_{(i,j)} > 0$ if $\lambda_* < 1$, where $\lambda_*$ is the largest eigenvalue of $(1+k)^{-2}T_hT_h'\tau T_h'T_h\tau^{-1}$. This completes the proof.

b) If $A \geq 0$, $\hat{y}_{SRrk}$ is superior to $\hat{y}_{MRE}$ in MSEM sense when the regression model is misspecified due to excluding relevant variables if and only if $\theta \in \mathfrak{R}(A)$ and $\theta' A^{-1}\theta \leq 1$, where $A = X_*\big(\sigma^2(\tau - (1+k)^{-2}T_hT_h'\tau T_h'T_h) + (\tau A)(\tau A)' - (1+k)^{-2}\big((T_hT_h' - (1+k)I)\gamma + T_hT_h'\tau A\big)\big((T_hT_h' - (1+k)I)\gamma + T_hT_h'\tau A\big)' + (1+k)^{-2}\big((1+k)I - T_hT_h'\big)(\gamma + \tau A)(\gamma + \tau A)'\big((1+k)I - T_hT_h'\big)'\big)X_*' + \delta\delta'$, $\mathfrak{R}(A)$ stands for column space of $A$ and $A^{-1}$ is an independent choice of g-inverse of $A$ and $\theta = \delta + (1+k)^{-1}X_*\big((1+k)I - T_hT_h'\big)(\gamma + \tau A)$.

**Corollary C7:**

a) If $\lambda_* < 1$, $\hat{\gamma}_{SRrd}$ is superior to $\hat{\gamma}_{MRE}$ in MSEM sense when the regression model is misspecified due to excluding relevant variables if and only if

$2^{-2}(1+d)^2\big((T_hT_h' - 2(1+d)^{-1}I)\gamma + T_hT_h'\tau A\big)'(\sigma^2(\tau - 2^{-2}(1+d)^2 T_hT_h'\tau T_h'T_h) + (\tau A)(\tau A)')^{-1}\big((T_hT_h' - 2(1+d)^{-1}I)\gamma + T_hT_h'\tau A\big) \leq 1$

where $\lambda_*$ is the largest eigenvalue of $2^{-2}(1+d)^2 T_hT_h'\tau T_h'T_h\tau^{-1}$.

**Proof:** Consider $D_{(i,j)} = D(\hat{\gamma}_{MRE}) - D(\hat{\gamma}_{SRrd}) = \sigma^2\tau - 2^{-2}(1+d)^2\sigma^2 T_hT_h'\tau T_h'T_h$
$= \sigma^2(\tau - 2^{-2}(1+d)^2 T_hT_h'\tau T_h'T_h)$

Since $\tau > 0$, according to Lemma A2 (see Appendix A) $D_{(i,j)} > 0$ if $\lambda_* < 1$, where $\lambda_*$ is the largest eigenvalue of $2^{-2}(1+d)^2 T_hT_h'\tau T_h'T_h\tau^{-1}$. This completes the proof.

b) If $A \geq 0$, $\hat{y}_{SRrd}$ is superior to $\hat{y}_{MRE}$ in MSEM sense when the regression model is misspecified due to excluding relevant variables if and only if $\theta \in \mathfrak{R}(A)$ and $\theta' A^{-1}\theta \leq 1$, where $A = X_*\big(\sigma^2(\tau - 2^{-2}(1+d)^2 T_hT_h'\tau T_h'T_h) + (\tau A)(\tau A)' - 2^{-2}(1+d)^2\big((T_hT_h' - 2(1+d)^{-1}I)\gamma + T_hT_h'\tau A\big)\big((T_hT_h' - 2(1+d)^{-1}I)\gamma + T_hT_h'\tau A\big)' + 2^{-2}(1+d)^2(2(1+d)^{-1}I - T_hT_h')(\gamma + \tau A)(\gamma + \tau A)'(2(1+d)^{-1}I - T_hT_h')'\big)X_*' + \delta\delta'$, $\mathfrak{R}(A)$ stands for column space of $A$ and $A^{-1}$ is an independent choice of g-inverse of $A$ and $\theta = \delta + 2^{-1}(1+d)X_*(2(1+d)^{-1}I - T_hT_h')(\gamma + \tau A)$.

**Corollary C8:**

a) $\hat{\gamma}_{SRRE}$ is superior to $\hat{\gamma}_{SRAURE}$ in MSEM sense when the regression model is misspecified due to excluding relevant variables if and only if

$(1+k)^2(\tau A - k\gamma)'\Big(k(2+3k)\sigma^2\tau + \big((1+2k)\tau A - k^2\gamma\big)\big((1+2k)\tau A - k^2\gamma\big)'\Big)^{-1}(\tau A - k\gamma) \leq 1$

**Proof:** Consider
$$D_{(i,j)} = D(\hat{\gamma}_{SRAURE}) - D(\hat{\gamma}_{SRRE}) = (1+k)^{-4}(1+2k)^2\sigma^2\tau - (1+k)^{-2}\sigma^2\tau$$
$$= ((1+2k)^2 - (1+k)^2)(1+k)^{-4}\sigma^2\tau$$
$$= k(2+3k)(1+k)^{-4}\sigma^2\tau$$

Since $k > 0$ and $\tau > 0$, hence $D_{(i,j)} > 0$. This completes the proof.

b) If $A \geq 0$, $\hat{y}_{SRRE}$ is superior to $\hat{y}_{SRAURE}$ in MSEM sense when the regression model is misspecified due to excluding relevant variables if and only if $\theta \in \Re(A)$ and $\theta'A^{-1}\theta \leq 1$, where $A = X_*\big(k(2+3k)(1+k)^{-4}\sigma^2\tau + (1+k)^{-4}\big((1+2k)\tau A - k^2\gamma\big)\big((1+2k)\tau A - k^2\gamma\big)' - (1+k)^{-2}(\tau A - k\gamma)(\tau A - k\gamma)' + k^2(1+k)^{-4}(\gamma + \tau A)(\gamma + \tau A)'\big)X'_* + \delta\delta'$, $\Re(A)$ stands for column space of $A$ and $A^{-1}$ is an independent choice of g-inverse of $A$ and $\theta = \delta + k(1+k)^{-2}X_*(\gamma + \tau A)$.

**Corollary C9:**

a) If $k > (1-d)(1+d)^{-1}$, $\hat{\gamma}_{SRRE}$ is superior to $\hat{\gamma}_{SRLE}$ in MSEM sense when the regression model is misspecified due to excluding relevant variables if and only if
$$2^2(1+k)^{-2}(\tau A - k\gamma)'\Big((1+k)^{-2}(k(1+d) + d + 3)(k(1+d) + d - 1)\sigma^2\tau + \big((1+d)\tau A - (1-d)\gamma\big)\big((1+d)\tau A - (1-d)\gamma\big)'\Big)^{-1}(\tau A - k\gamma) \leq 1$$

**Proof:** Consider $D_{(i,j)} = D(\hat{\gamma}_{SRLE}) - D(\hat{\gamma}_{SRRE}) = 2^{-2}(1+d)^2\sigma^2\tau - (1+k)^{-2}\sigma^2\tau$
$$= ((1+k)^2(1+d)^2 - 2^2)2^{-2}(1+k)^{-2}\sigma^2\tau$$
$$= 2^{-2}(1+k)^{-2}(k(1+d) + d + 3)(k(1+d) + d - 1)\sigma^2\tau$$

Since $k > 0$, $0 < d < 1$ and $\tau > 0$. $D_{(i,j)} > 0$ if $(k(1+d) + d - 1) > 0$, which implies $k > (1-d)(1+d)^{-1}$. This complete the proof.

b) If $A \geq 0$, $\hat{y}_{SRRE}$ is superior to $\hat{y}_{SRLE}$ in MSEM sense when the regression model is misspecified due to excluding relevant variables if and only if $\theta \in \Re(A)$ and $\theta'A^{-1}\theta \leq 1$, where $A = X_*\big(2^{-2}(1+k)^{-2}(k(1+d)+d+3)(k(1+d)+d-1)\sigma^2\tau + 2^{-2}\big((1+d)\tau A - (1-d)\gamma\big)\big((1+d)\tau A - (1-d)\gamma\big)' - (1+k)^{-2}(\tau A - k\gamma)(\tau A - k\gamma)' + 2^{-2}(1+k)^{-2}(1+k+d+kd)^2(\gamma + \tau A)(\gamma + \tau A)'\big)X'_* + \delta\delta'$, $\Re(A)$ stands for column space of $A$ and $A^{-1}$ is an independent choice of g-inverse of $A$ and $\theta = \delta + 2^{-1}(1+k)^{-1}(1+k+d+kd)X_*(\gamma + \tau A)$.

**Corollary C10:**

a) If $k > (1-d)^2(1+d)^{-1}(3-d)^{-1}$, $\hat{\gamma}_{SRRE}$ is superior to $\hat{\gamma}_{SRAULE}$ in MSEM sense when the regression model is misspecified due to excluding relevant variables if and only if
$$2^4(1+k)^{-2}(\tau A - k\gamma)'\Big((1+k)^{-2}(7+2d-d^2+k(1+d)(3-d))(k(1+d)(3-d) - (1-d)^2)\sigma^2\tau + \big((1+d)(3-d)\tau A - (1-d)^2\gamma\big)\big((1+d)(3-d)\tau A - (1-d)^2\gamma\big)'\Big)^{-1}(\tau A - k\gamma) \leq 1$$

**Proof:** Consider

$$D_{(i,j)} = D(\hat{\gamma}_{SRAULE}) - D(\hat{\gamma}_{SRRE}) = 2^{-4}(1+d)^2(3-d)^2\sigma^2\tau - (1+k)^{-2}\sigma^2\tau$$
$$= ((1+d)^2(3-d)^2(1+k)^2 - 2^4)2^{-4}(1+k)^{-2}\sigma^2\tau$$
$$= 2^{-4}(1+k)^{-2}(7+2d-d^2+k(1+d)(3-d))(k(1+d)(3-d)-(1-d)^2)\sigma^2\tau$$

Since $k > 0$, $0 < d < 1$ and $\tau > 0$. $D_{(i,j)} > 0$ if $(k(1+d)(3-d) - (1-d)^2) > 0$, which implies $k > (1-d)^2(1+d)^{-1}(3-d)^{-1}$. This completes the proof.

b) If $A \geq 0$, $\hat{y}_{SRRE}$ is superior to $\hat{y}_{SRAULE}$ in MSEM sense when the regression model is misspecified due to excluding relevant variables if and only if $\theta \in \Re(A)$ and $\theta' A^{-1}\theta \leq 1$, where $A = X_* \big(2^{-4}(1+k)^{-2}(7+2d-d^2+k(1+d)(3-d))(k(1+d)(3-d)-(1-d)^2)\sigma^2\tau + 2^{-4}((1+d)(3-d)\tau A - (1-d)^2\gamma)((1+d)(3-d)\tau A - (1-d)^2\gamma)' - (1+k)^{-2}(\tau A - k\gamma)(\tau A - k\gamma)' + 2^{-4}(1+k)^{-2}((1+k)d(2-d)+3k-1)^2(\gamma+\tau A)(\gamma+\tau A)'\big)X_*' + \delta\delta'$, $\Re(A)$ stands for column space of $A$ and $A^{-1}$ is an independent choice of g-inverse of $A$ and $\theta = \delta + 2^{-2}(1+k)^{-1}((1+k)d(2-d)+3k-1)X_*(\gamma+\tau A)$.

**Corollary C11:**
a) If $\lambda_* < 1$, $\hat{\gamma}_{SRPCR}$ is superior to $\hat{\gamma}_{SRRE}$ in MSEM sense when the regression model is misspecified due to excluding relevant variables if and only if
$$\big((T_h T_h' - I)\gamma + T_h T_h'\tau A\big)' \big(((1+k)^{-2}\tau - T_h T_h'\tau T_h' T_h)\sigma^2 + (1+k)^{-2}(\tau A - k\gamma)(\tau A - k\gamma)'\big)^{-1}\big((T_h T_h' - I)\gamma + T_h T_h'\tau A\big) \leq 1$$
where $\lambda_*$ is the largest eigenvalue of $(1+k)^2 T_h T_h'\tau T_h' T_h \tau^{-1}$.

**Proof:** Consider $D_{(i,j)} = D(\hat{\gamma}_{SRRE}) - D(\hat{\gamma}_{SRPCR}) = (1+k)^{-2}\sigma^2\tau - \sigma^2 T_h T_h'\tau T_h' T_h$
$$= ((1+k)^{-2}\tau - T_h T_h'\tau T_h' T_h)\sigma^2$$

Since $\tau > 0$, according to Lemma A2 (see Appendix A) $D_{(i,j)} > 0$ if $\lambda_* < 1$, where $\lambda_*$ is the largest eigenvalue of $(1+k)^2 T_h T_h'\tau T_h' T_h \tau^{-1}$. This completes the proof.

b) If $A \geq 0$, $\hat{y}_{SRPCR}$ is superior to $\hat{y}_{SRRE}$ in MSEM sense when the regression model is misspecified due to excluding relevant variables if and only if $\theta \in \Re(A)$ and $\theta' A^{-1}\theta \leq 1$, where $A = X_* \big(((1+k)^{-2}\tau - T_h T_h'\tau T_h' T_h)\sigma^2 + (1+k)^{-2}(\tau A - k\gamma)(\tau A - k\gamma)' - \big((T_h T_h' - I)\gamma + T_h T_h'\tau A\big)\big((T_h T_h' - I)\gamma + T_h T_h'\tau A\big)' + (1+k)^{-2}(I - (1+k)T_h T_h')(\gamma+\tau A)(\gamma+\tau A)'(I - (1+k)T_h T_h')'\big)X_*' + \delta\delta'$, $\Re(A)$ stands for column space of $A$ and $A^{-1}$ is an independent choice of g-inverse of $A$ and $\theta = \delta + (1+k)^{-1}X_*(I - (1+k)T_h T_h')(\gamma+\tau A)$.

**Corollary C12:**
a) If $\lambda_* < 1$, $\hat{\gamma}_{SRrk}$ is superior to $\hat{\gamma}_{SRRE}$ in MSEM sense when the regression model is misspecified due to excluding relevant variables if and only if
$$\big((T_h T_h' - (1+k)I)\gamma + T_h T_h'\tau A\big)' \big((\tau - T_h T_h'\tau T_h' T_h)\sigma^2 + (\tau A - k\gamma)(\tau A - k\gamma)'\big)^{-1}\big((T_h T_h' - (1+k)I)\gamma + T_h T_h'\tau A\big) \leq 1$$
where $\lambda_*$ is the largest eigenvalue of $T_h T_h'\tau T_h' T_h \tau^{-1}$.

**Proof:** Consider $D_{(i,j)} = D(\hat{y}_{SRRE}) - D(\hat{y}_{SRrk}) = (1+k)^{-2}\sigma^2\tau - (1+k)^{-2}\sigma^2 T_h T_h' \tau T_h' T_h$
$$= (\tau - T_h T_h' \tau T_h' T_h)(1+k)^{-2}\sigma^2$$

Since $\tau > 0$, according to Lemma A2 (see Appendix A) $D_{(i,j)} > 0$ if $\lambda_* < 1$, where $\lambda_*$ is the largest eigenvalue of $T_h T_h' \tau T_h' T_h \tau^{-1}$. This completes the proof.

b) If $A \geq 0$, $\hat{y}_{SRrk}$ is superior to $\hat{y}_{SRRE}$ in MSEM sense when the regression model is misspecified due to excluding relevant variables if and only if $\theta \in \mathfrak{R}(A)$ and $\theta'A^{-1}\theta \leq 1$, where $A = X_*\Big((\tau - T_h T_h' \tau T_h' T_h)(1+k)^{-2}\sigma^2 + (1+k)^{-2}(\tau A - k\gamma)(\tau A - k\gamma)' - (1+k)^{-2}\big((T_h T_h' - (1+k)I)\gamma + T_h T_h' \tau A\big)\big((T_h T_h' - (1+k)I)\gamma + T_h T_h' \tau A\big)' + (1+k)^{-2}(I - T_h T_h')(\gamma + \tau A)(\gamma + \tau A)'(I - T_h T_h')'\Big)X_*' + \delta\delta'$, $\mathfrak{R}(A)$ stands for column space of $A$ and $A^{-1}$ is an independent choice of g-inverse of $A$ and $\theta = \delta + (1+k)^{-1}X_*(I - T_h T_h')(\gamma + \tau A)$.

**Corollary C13:**

a) If $\lambda_* < 1$, $\hat{\gamma}_{SRrd}$ is superior to $\hat{\gamma}_{SRRE}$ in MSEM sense when the regression model is misspecified due to excluding relevant variables if and only if
$2^{-2}(1+d)^2\big((T_h T_h' - 2(1+d)^{-1}I)\gamma + T_h T_h' \tau A\big)'\big(((1+k)^{-2}\tau - 2^{-2}(1+d)^2 T_h T_h' \tau T_h' T_h)\sigma^2 + (1+k)^{-2}(\tau A - k\gamma)(\tau A - k\gamma)'\big)^{-1}\big((T_h T_h' - 2(1+d)^{-1}I)\gamma + T_h T_h' \tau A\big) \leq 1$
where $\lambda_*$ is the largest eigenvalue of $2^{-2}(1+d)^2(1+k)^2 T_h T_h' \tau T_h' T_h \tau^{-1}$.

**Proof:** Consider
$$D_{(i,j)} = D(\hat{y}_{SRRE}) - D(\hat{y}_{SRrd}) = (1+k)^{-2}\sigma^2\tau - 2^{-2}(1+d)^2\sigma^2 T_h T_h' \tau T_h' T_h$$
$$= ((1+k)^{-2}\tau - 2^{-2}(1+d)^2 T_h T_h' \tau T_h' T_h)\sigma^2$$

Since $\tau > 0$, according to Lemma A2 (see Appendix A) $D_{(i,j)} > 0$ if $\lambda_* < 1$, where $\lambda_*$ is the largest eigenvalue of $2^{-2}(1+d)^2(1+k)^2 T_h T_h' \tau T_h' T_h \tau^{-1}$. This completes the proof.

b) If $A \geq 0$, $\hat{y}_{SRrd}$ is superior to $\hat{y}_{SRRE}$ in MSEM sense when the regression model is misspecified due to excluding relevant variables if and only if $\theta \in \mathfrak{R}(A)$ and $\theta'A^{-1}\theta \leq 1$, where $A = X_*\Big(((1+k)^{-2}\tau - 2^{-2}(1+d)^2 T_h T_h' \tau T_h' T_h)\sigma^2 + (1+k)^{-2}(\tau A - k\gamma)(\tau A - k\gamma)' - 2^{-2}(1+d)^2\big((T_h T_h' - 2(1+d)^{-1}I)\gamma + T_h T_h' \tau A\big)\big((T_h T_h' - 2(1+d)^{-1}I)\gamma + T_h T_h' \tau A\big)' + 2^{-2}(1+k)^{-2}(2I - (1+k)(1+d)T_h T_h')(\gamma + \tau A)(\gamma + \tau A)'(2I - (1+k)(1+d)T_h T_h')'\Big)X_*' + \delta\delta'$, $\mathfrak{R}(A)$ stands for column space of $A$ and $A^{-1}$ is an independent choice of g-inverse of $A$ and $\theta = \delta + 2^{-1}(1+k)^{-1}X_*(2I - (1+k)(1+d)T_h T_h')(\gamma + \tau A)$.

**Corollary C14:**

a) If $d > (1 + 2k - k^2)(1+k)^{-2}$, $\hat{\gamma}_{SRAURE}$ is superior to $\hat{\gamma}_{SRLE}$ in MSEM sense when the regression model is misspecified due to excluding relevant variables if and only if

$$2^2\left((1+2k)\tau A - k^2\gamma\right)'\Big((3+6k+k^2+d(1+k)^2)(k^2-2k-1+d(1+k)^2)\sigma^2\tau +$$
$$(1+k)^4\big((1+d)\tau A - (1-d)\gamma\big)\big((1+d)\tau A - (1-d)\gamma\big)'\Big)^{-1}\left((1+2k)\tau A - k^2\gamma\right) \leq 1$$

**Proof:** Consider
$$D_{(i,j)} = D(\hat{y}_{SRLE}) - D(\hat{y}_{SRAURE}) = 2^{-2}(1+d)^2\sigma^2\tau - (1+k)^{-4}(1+2k)^2\sigma^2\tau$$
$$= \big((1+d)^2(1+k)^4 - 2^2(1+2k)^2\big)2^{-2}(1+k)^{-4}\sigma^2\tau$$
$$= (3+6k+k^2+d(1+k)^2)(k^2-2k-1+d(1+k)^2)2^{-2}(1+k)^{-4}\sigma^2\tau$$

Since $k > 0$, $0 < d < 1$ and $\tau > 0$. $D_{(i,j)} > 0$ if $(k^2 - 2k - 1 + d(1+k)^2) > 0$, which implies $d > (1 + 2k - k^2)(1+k)^{-2}$. This completes the proof.

b) If $A \geq 0$, $\hat{y}_{SRAURE}$ is superior to $\hat{y}_{SRLE}$ in MSEM sense when the regression model is misspecified due to excluding relevant variables if and only if $\theta \in \mathfrak{R}(A)$ and $\theta'A^{-1}\theta \leq 1$, where $A = X_*\Big((3+6k+k^2+d(1+k)^2)(k^2-2k-1+d(1+k)^2)2^{-2}(1+k)^{-4}\sigma^2\tau + 2^{-2}\big((1+d)\tau A - (1-d)\gamma\big)\big((1+d)\tau A - (1-d)\gamma\big)' - (1+k)^{-4}\big((1+2k)\tau A - k^2\gamma\big)\big((1+2k)\tau A - k^2\gamma\big)' + 2^{-2}(1+k)^{-4}\big((1+d)k^2 + (2k+1)(d-1)\big)^2(\gamma+\tau A)(\gamma+\tau A)'\Big)X'_* + \delta\delta'$, $\mathfrak{R}(A)$ stands for column space of $A$ and $A^{-1}$ is an independent choice of g-inverse of $A$ and $\theta = \delta + 2^{-1}(1+k)^{-2}\big((1+d)k^2 + (2k+1)(d-1)\big)X_*(\gamma + \tau A)$.

**Corollary C15:**
a) If $(1+d)(3-d) > 4(1+2k)(1+k)^{-2}$, $\hat{\gamma}_{SRAURE}$ is superior to $\hat{\gamma}_{SRAULE}$ in MSEM sense when the regression model is misspecified due to excluding relevant variables if and only if
$$2^4\left((1+2k)\tau A - k^2\gamma\right)'\Big(\big((1+k)^2(1+d)(3-d) + 4(1+2k)\big)\big((1+k)^2(1+d)(3-d) -$$
$$4(1+2k)\big)\sigma^2\tau + (1+k)^4\big((1+d)(3-d)\tau A - (1-d)^2\gamma\big)\big((1+d)(3-d)\tau A - (1-$$
$$d)^2\gamma\big)'\Big)^{-1}\left((1+2k)\tau A - k^2\gamma\right) \leq 1$$

**Proof:** Consider
$$D_{(i,j)} = D(\hat{\gamma}_{SRAULE}) - D(\hat{\gamma}_{SRAURE}) = 2^{-4}(1+d)^2(3-d)^2\sigma^2\tau - (1+k)^{-4}(1+2k)^2\sigma^2\tau$$
$$= \big((1+k)^4(1+d)^2(3-d)^2 - 2^4(1+2k)^2\big)2^{-4}(1+k)^{-4}\sigma^2\tau$$
$$= \big((1+k)^2(1+d)(3-d) + 4(1+2k)\big)\big((1+k)^2(1+d)(3-d) - 4(1+2k)\big)$$
$$2^{-4}(1+k)^{-4}\sigma^2\tau$$

Since $k > 0$, $0 < d < 1$ and $\tau > 0$. $D_{(i,j)} > 0$ if $\big((1+k)^2(1+d)(3-d) - 4(1+2k)\big) > 0$, which implies $(1+d)(3-d) > 4(1+2k)(1+k)^{-2}$. This completes the proof.

b) If $A \geq 0$, $\hat{y}_{SRAURE}$ is superior to $\hat{y}_{SRAULE}$ in MSEM sense when the regression model is misspecified due to excluding relevant variables if and only if $\theta \in \mathfrak{R}(A)$ and $\theta'A^{-1}\theta \leq 1$, where $A = X_*\Big(\big((1+k)^2(1+d)(3-d) + 4(1+2k)\big)\big((1+k)^2(1+d)(3-d) - 4(1+2k)\big)2^{-4}(1+k)^{-4}\sigma^2\tau + 2^{-4}\big((1+d)(3-d)\tau A - (1-d)^2\gamma\big)\big((1+d)(3-d)\tau A - (1-d)^2\gamma\big)' - (1+$ $k)^{-4}\left((1+2k)\tau A-k^2\gamma\right)\left((1+2k)\tau A-k^2\gamma\right)'+2^{-4}(1+k)^{-4}\big((k-1)(3k+1)+d(2-d)(1+k)^2\big)^2(\gamma+\tau A)(\gamma+\tau A)'\Big)X'_*+\delta\delta'$, $\mathfrak{R}(A)$ stands for column space of $A$ and $A^{-1}$ is an independent choice of g-inverse of $A$ and $\theta=\delta+2^{-2}(1+k)^{-2}\big((k-1)(3k+1)+d(2-d)(1+k)^2\big)X_*(\gamma+\tau A)$.

**Corollary C16:**

a) If $\lambda_*<1$, $\hat{\gamma}_{SRPCR}$ is superior to $\hat{\gamma}_{SRAURE}$ in MSEM sense when the regression model is misspecified due to excluding relevant variables if and only if

$$\big((T_hT'_h-I)\gamma+T_hT'_h\tau A\big)'\Big(((1+k)^{-4}(1+2k)^2\tau-T_hT'_h\tau T'_hT_h)\sigma^2+(1+k)^{-4}\big((1+2k)\tau A-k^2\gamma\big)\big((1+2k)\tau A-k^2\gamma\big)'\Big)^{-1}\big((T_hT'_h-I)\gamma+T_hT'_h\tau A\big)\leq 1$$

where $\lambda_*$ is the largest eigenvalue of $(1+k)^4(1+2k)^{-2}T_hT'_h\tau T'_hT_h\tau^{-1}$.

**Proof:** Consider

$$D_{(i,j)}=D(\hat{\gamma}_{SRAURE})-D(\hat{\gamma}_{SRPCR})=(1+k)^{-4}(1+2k)^2\sigma^2\tau-\sigma^2T_hT'_h\tau T'_hT_h$$
$$=((1+k)^{-4}(1+2k)^2\tau-T_hT'_h\tau T'_hT_h)\sigma^2$$

Since $\tau>0$, according to Lemma A2 (see Appendix A) $D_{(i,j)}>0$ if $\lambda_*<1$, where $\lambda_*$ is the largest eigenvalue of $(1+k)^4(1+2k)^{-2}T_hT'_h\tau T'_hT_h\tau^{-1}$. This completes the proof.

b) If $A\geq 0$, $\hat{\gamma}_{SRPCR}$ is superior to $\hat{\gamma}_{SRAURE}$ in MSEM sense when the regression model is misspecified due to excluding relevant variables if and only if $\theta\in\mathfrak{R}(A)$ and $\theta'A^{-1}\theta\leq 1$, where $A=X_*\Big(((1+k)^{-4}(1+2k)^2\tau-T_hT'_h\tau T'_hT_h)\sigma^2+(1+k)^{-4}\big((1+2k)\tau A-k^2\gamma\big)\big((1+2k)\tau A-k^2\gamma\big)'-\big((T_hT'_h-I)\gamma+T_hT'_h\tau A\big)\big((T_hT'_h-I)\gamma+T_hT'_h\tau A\big)'+(1+k)^{-4}\big((1+2k)I-(1+k)^2T_hT'_h\big)(\gamma+\tau A)(\gamma+\tau A)'\big((1+2k)I-(1+k)^2T_hT'_h\big)'\Big)X'_*+\delta\delta'$, $\mathfrak{R}(A)$ stands for column space of $A$ and $A^{-1}$ is an independent choice of g-inverse of $A$ and $\theta=\delta+(1+k)^{-2}X_*\big((1+2k)I-(1+k)^2T_hT'_h\big)(\gamma+\tau A)$.

**Corollary C17:**

a) If $\lambda_*<1$, $\hat{\gamma}_{SRrk}$ is superior to $\hat{\gamma}_{SRAURE}$ in MSEM sense when the regression model is misspecified due to excluding relevant variables if and only if

$$(1+k)^2\big((T_hT'_h-(1+k)I)\gamma+T_hT'_h\tau A\big)'\Big(((1+2k)^2\tau-(1+k)^2T_hT'_h\tau T'_hT_h)\sigma^2+\big((1+2k)\tau A-k^2\gamma\big)\big((1+2k)\tau A-k^2\gamma\big)'\Big)^{-1}\big((T_hT'_h-(1+k)I)\gamma+T_hT'_h\tau A\big)\leq 1$$

where $\lambda_*$ is the largest eigenvalue of $(1+2k)^{-2}(1+k)^2T_hT'_h\tau T'_hT_h\tau^{-1}$.

**Proof:** Consider

$$D_{(i,j)}=D(\hat{\gamma}_{SRAURE})-D(\hat{\gamma}_{SRrk})=(1+k)^{-4}(1+2k)^2\sigma^2\tau-(1+k)^{-2}\sigma^2T_hT'_h\tau T'_hT_h$$
$$=((1+2k)^2\tau-(1+k)^2T_hT'_h\tau T'_hT_h)(1+k)^{-4}\sigma^2$$

Since $\tau > 0$, according to Lemma A2 (see Appendix A) $D_{(i,j)} > 0$ if $\lambda_* < 1$, where $\lambda_*$ is the largest eigenvalue of $(1 + 2k)^{-2}(1 + k)^2 T_h T_h' \tau T_h' T_h \tau^{-1}$. This completes the proof.

b) If $A \geq 0$, $\hat{y}_{SRrk}$ is superior to $\hat{y}_{SRAURE}$ in MSEM sense when the regression model is misspecified due to excluding relevant variables if and only if $\theta \in \mathfrak{R}(A)$ and $\theta' A^{-1}\theta \leq 1$, where $A = X_*\Big(((1 + 2k)^2\tau - (1 + k)^2 T_h T_h' \tau T_h' T_h)(1 + k)^{-4}\sigma^2 + (1 + k)^{-4}\big((1 + 2k)\tau A - k^2\gamma\big)\big((1 + 2k)\tau A - k^2\gamma\big)' - (1 + k)^{-2}\big((T_h T_h' - (1 + k)I)\gamma + T_h T_h' \tau A\big)\big((T_h T_h' - (1 + k)I)\gamma + T_h T_h' \tau A\big)' + (1 + k)^{-4}\big((1 + 2k)I - (1 + k)T_h T_h'\big)(\gamma + \tau A)(\gamma + \tau A)'\big((1 + 2k)I - (1 + k)T_h T_h'\big)'\Big)X_*' + \delta\delta'$,

$\mathfrak{R}(A)$ stands for column space of $A$ and $A^{-1}$ is an independent choice of g-inverse of $A$ and $\theta = \delta + (1 + k)^{-2} X_*\big((1 + 2k)I - (1 + k)T_h T_h'\big)(\gamma + \tau A)$.

**Corollary C18:**

a) If $\lambda_* < 1$, $\hat{\gamma}_{SRrd}$ is superior to $\hat{\gamma}_{SRAURE}$ in MSEM sense when the regression model is misspecified due to excluding relevant variables if and only if

$2^{-2}(1 + d)^2\big((T_h T_h' - 2(1 + d)^{-1}I)\gamma + T_h T_h' \tau A\big)' \Big(((1 + k)^{-4}(1 + 2k)^2\tau - 2^{-2}(1 + d)^2 T_h T_h' \tau T_h' T_h)\sigma^2 + (1 + k)^{-4}\big((1 + 2k)\tau A - k^2\gamma\big)\big((1 + 2k)\tau A - k^2\gamma\big)'\Big)^{-1} \big((T_h T_h' - 2(1 + d)^{-1}I)\gamma + T_h T_h' \tau A\big) \leq 1$

where $\lambda_*$ is the largest eigenvalue of $2^{-2}(1 + d)^2(1 + k)^4(1 + 2k)^{-2}T_h T_h' \tau T_h' T_h \tau^{-1}$.

**Proof:** Consider
$D_{(i,j)} = D(\hat{\gamma}_{SRAURE}) - D(\hat{\gamma}_{SRrk}) = (1 + k)^{-4}(1 + 2k)^2\sigma^2\tau - 2^{-2}(1 + d)^2\sigma^2 T_h T_h' \tau T_h' T_h$
$= \big((1 + k)^{-4}(1 + 2k)^2\tau - 2^{-2}(1 + d)^2 T_h T_h' \tau T_h' T_h\big)\sigma^2$

Since $\tau > 0$, according to Lemma A2 (see Appendix A) $D_{(i,j)} > 0$ if $\lambda_* < 1$, where $\lambda_*$ is the largest eigenvalue of $2^{-2}(1 + d)^2(1 + k)^4(1 + 2k)^{-2}T_h T_h' \tau T_h' T_h \tau^{-1}$. This completes the proof.

b) If $A \geq 0$, $\hat{y}_{SRrd}$ is superior to $\hat{y}_{SRAURE}$ in MSEM sense when the regression model is misspecified due to excluding relevant variables if and only if $\theta \in \mathfrak{R}(A)$ and $\theta' A^{-1}\theta \leq 1$, where $A = X_*\Big(((1 + k)^{-4}(1 + 2k)^2\tau - 2^{-2}(1 + d)^2 T_h T_h' \tau T_h' T_h)\sigma^2 + (1 + k)^{-4}\big((1 + 2k)\tau A - k^2\gamma\big)\big((1 + 2k)\tau A - k^2\gamma\big)' - 2^{-2}(1 + d)^2\big((T_h T_h' - 2(1 + d)^{-1}I)\gamma + T_h T_h' \tau A\big)\big((T_h T_h' - 2(1 + d)^{-1}I)\gamma + T_h T_h' \tau A\big)' + 2^{-2}(1 + k)^{-4}\big(2(1 + 2k)I - (1 + k)^2(1 + d)T_h T_h'\big)(\gamma + \tau A)(\gamma + \tau A)'\big(2(1 + 2k)I - (1 + k)^2(1 + d)T_h T_h'\big)'\Big)X_*' + \delta\delta'$, $\mathfrak{R}(A)$ stands for column space of $A$ and $A^{-1}$ is an independent choice of g-inverse of $A$ and $\theta = \delta + 2^{-1}(1 + k)^{-2}X_*\big(2(1 + 2k)I - (1 + k)^2(1 + d)T_h T_h'\big)(\gamma + \tau A)$.

**Corollary C19:**

a) $\hat{\gamma}_{SRLE}$ is superior to $\hat{\gamma}_{SRAULE}$ in MSEM sense when the regression model is misspecified due to excluding relevant variables if and only if

$2^2\big((1 + d)\tau A - (1 - d)\gamma\big)' \Big((5 - d)(1 - d)(1 + d)^2\sigma^2\tau + \big((1 + d)(3 - d)\tau A - (1 - d)^2\gamma\big)\big((1 + d)(3 - d)\tau A - (1 - d)^2\gamma\big)'\Big)^{-1} \big((1 + d)\tau A - (1 - d)\gamma\big) \leq 1$

**Proof:** Consider $D_{(i,j)} = D(\hat{y}_{SRAULE}) - D(\hat{y}_{SRLE}) = 2^{-4}(1+d)^2(3-d)^2\sigma^2\tau - 2^{-2}(1+d)^2\sigma^2\tau$
$$= ((3-d)^2 - 2^2)2^{-4}(1+d)^2\sigma^2\tau$$
$$= (5-d)(1-d)2^{-4}(1+d)^2\sigma^2\tau$$

Since $0 < d < 1$ and $\tau > 0$, hence $D_{(i,j)} > 0$. This completes the proof.

b) If $A \geq 0$, $\hat{y}_{SRLE}$ is superior to $\hat{y}_{SRAULE}$ in MSEM sense when the regression model is misspecified due to excluding relevant variables if and only if $\theta \in \Re(A)$ and $\theta'A^{-1}\theta \leq 1$, where $A = X_*\Big((5-d)(1-d)2^{-4}(1+d)^2\sigma^2\tau + 2^{-4}((1+d)(3-d)\tau A - (1-d)^2\gamma)((1+d)(3-d)\tau A - (1-d)^2\gamma)' - 2^{-2}((1+d)\tau A - (1-d)\gamma)((1+d)\tau A - (1-d)\gamma)' + 2^{-4}(1+d)^2(1-d)^2(\gamma + \tau A)(\gamma + \tau A)'\Big)X_*' + \delta\delta'$, $\Re(A)$ stands for column space of $A$ and $A^{-1}$ is an independent choice of g-inverse of $A$ and $\theta = \delta + 2^{-2}(1+d)(1-d)X_*(\gamma + \tau A)$.

**Corollary C20:**

a) If $\lambda_* < 1$, $\hat{\gamma}_{SRPCR}$ is superior to $\hat{\gamma}_{SRLE}$ in MSEM sense when the regression model is misspecified due to excluding relevant variables if and only if

$$((T_hT_h' - I)\gamma + T_hT_h'\tau A)'\Big((2^{-2}(1+d)^2\tau - T_hT_h'\tau T_h'T_h)\sigma^2 + 2^{-2}((1+d)\tau A - (1-d)\gamma)((1+d)\tau A - (1-d)\gamma)'\Big)^{-1}\big((T_hT_h' - I)\gamma + T_hT_h'\tau A\big) \leq 1$$

where $\lambda_*$ is the largest eigenvalue of $2^2(1+d)^{-2}T_hT_h'\tau T_h'T_h\tau^{-1}$.

**Proof:** Consider $D_{(i,j)} = D(\hat{y}_{SRLE}) - D(\hat{y}_{SRPCR}) = 2^{-2}(1+d)^2\sigma^2\tau - \sigma^2 T_hT_h'\tau T_h'T_h$
$$= (2^{-2}(1+d)^2\tau - T_hT_h'\tau T_h'T_h)\sigma^2$$

Since $\tau > 0$, according to Lemma A2 (see Appendix A) $D_{(i,j)} > 0$ if $\lambda_* < 1$, where $\lambda_*$ is the largest eigenvalue of $2^2(1+d)^{-2}T_hT_h'\tau T_h'T_h\tau^{-1}$. This completes the proof.

b) If $A \geq 0$, $\hat{y}_{SRPCR}$ is superior to $\hat{y}_{SRLE}$ in MSEM sense when the regression model is misspecified due to excluding relevant variables if and only if $\theta \in \Re(A)$ and $\theta'A^{-1}\theta \leq 1$, where $A = X_*\Big((2^{-2}(1+d)^2\tau - T_hT_h'\tau T_h'T_h)\sigma^2 + 2^{-2}((1+d)\tau A - (1-d)\gamma)((1+d)\tau A - (1-d)\gamma)' - ((T_hT_h' - I)\gamma + T_hT_h'\tau A)((T_hT_h' - I)\gamma + T_hT_h'\tau A)' + 2^{-2}((1+d)I - 2T_hT_h')(\gamma + \tau A)(\gamma + \tau A)'((1+d)I - 2T_hT_h')'\Big)X_*' + \delta\delta'$, $\Re(A)$ stands for column space of $A$ and $A^{-1}$ is an independent choice of g-inverse of $A$ and $\theta = \delta + 2^{-1}X_*((1+d)I - 2T_hT_h')(\gamma + \tau A)$.

**Corollary C21:**

a) If $\lambda_* < 1$, $\hat{\gamma}_{SRrk}$ is superior to $\hat{\gamma}_{SRLE}$ in MSEM sense when the regression model is misspecified due to excluding relevant variables if and only if

$$(1+k)^{-2}\big((T_hT_h' - (1+k)I)\gamma + T_hT_h'\tau A\big)'\Big((2^{-2}(1+d)^2\tau - (1+k)^{-2}T_hT_h'\tau T_h'T_h)\sigma^2 + 2^{-2}((1+d)\tau A - (1-d)\gamma)((1+d)\tau A - (1-d)\gamma)'\Big)^{-1}\big((T_hT_h' - (1+k)I)\gamma + T_hT_h'\tau A\big) \leq 1$$

where $\lambda_*$ is the largest eigenvalue of $2^2(1+d)^{-2}(1+k)^{-2}T_hT_h'\tau T_h'T_h\tau^{-1}$.

**Proof:** Consider $D_{(i,j)} = D(\hat{y}_{SRLE}) - D(\hat{y}_{SRrk}) = 2^{-2}(1+d)^2\sigma^2\tau - (1+k)^{-2}\sigma^2 T_h T_h' \tau T_h' T_h$
$$= (2^{-2}(1+d)^2\tau - (1+k)^{-2} T_h T_h' \tau T_h' T_h)\sigma^2$$

Since $\tau > 0$, according to Lemma A2 (see Appendix A) $D_{(i,j)} > 0$ if $\lambda_* < 1$, where $\lambda_*$ is the largest eigenvalue of $2^2(1+d)^{-2}(1+k)^{-2} T_h T_h' \tau T_h' T_h \tau^{-1}$. This completes the proof.

b) If $A \geq 0$, $\hat{y}_{SRrk}$ is superior to $\hat{y}_{SRLE}$ in MSEM sense when the regression model is misspecified due to excluding relevant variables if and only if $\theta \in \mathfrak{R}(A)$ and $\theta' A^{-1} \theta \leq 1$, where $A = X_*\Big((2^{-2}(1+d)^2\tau - (1+k)^{-2} T_h T_h' \tau T_h' T_h)\sigma^2 + 2^{-2}\big((1+d)\tau A - (1-d)\gamma\big)\big((1+d)\tau A - (1-d)\gamma\big)' - (1+k)^{-2}\big((T_h T_h' - (1+k)I)\gamma + T_h T_h' \tau A\big)\big((T_h T_h' - (1+k)I)\gamma + T_h T_h' \tau A\big)' + 2^{-2}(1+k)^{-2}\big((1+k)(1+d)I - 2T_h T_h'\big)(\gamma + \tau A)(\gamma + \tau A)'\big((1+k)(1+d)I - 2T_h T_h'\big)'\Big)X_*' + \delta\delta'$, $\mathfrak{R}(A)$ stands for column space of $A$ and $A^{-1}$ is an independent choice of g-inverse of $A$ and $\theta = \delta + 2^{-1}(1+k)^{-1}X_*((1+k)(1+d)I - 2T_h T_h')(\gamma + \tau A)$.

**Corollary C22:**

a) If $\lambda_* < 1$, $\hat{y}_{SRrd}$ is superior to $\hat{y}_{SRLE}$ in MSEM sense when the regression model is misspecified due to excluding relevant variables if and only if

$$\big((T_h T_h' - 2(1+d)^{-1}I)\gamma + T_h T_h' \tau A\big)'\Big((\tau - T_h T_h' \tau T_h' T_h)\sigma^2 + (1+d)^{-2}\big((1+d)\tau A - (1-d)\gamma\big)\big((1+d)\tau A - (1-d)\gamma\big)'\Big)^{-1}\big((T_h T_h' - 2(1+d)^{-1}I)\gamma + T_h T_h' \tau A\big) \leq 1$$

where $\lambda_*$ is the largest eigenvalue of $T_h T_h' \tau T_h' T_h \tau^{-1}$.

**Proof:** Consider $D_{(i,j)} = D(\hat{y}_{SRLE}) - D(\hat{y}_{SRrd}) = 2^{-2}(1+d)^2\sigma^2\tau - 2^{-2}(1+d)^2\sigma^2 T_h T_h' \tau T_h' T_h$
$$= (\tau - T_h T_h' \tau T_h' T_h) 2^{-2}(1+d)^2 \sigma^2$$

Since $\tau > 0$, according to Lemma A2 (see Appendix A) $D_{(i,j)} > 0$ if $\lambda_* < 1$, where $\lambda_*$ is the largest eigenvalue of $T_h T_h' \tau T_h' T_h \tau^{-1}$. This completes the proof.

b) If $A \geq 0$, $\hat{y}_{SRrd}$ is superior to $\hat{y}_{SRLE}$ in MSEM sense when the regression model is misspecified due to excluding relevant variables if and only if $\theta \in \mathfrak{R}(A)$ and $\theta' A^{-1} \theta \leq 1$, where $A = X_*\Big((\tau - T_h T_h' \tau T_h' T_h)2^{-2}(1+d)^2\sigma^2 + 2^{-2}\big((1+d)\tau A - (1-d)\gamma\big)\big((1+d)\tau A - (1-d)\gamma\big)' - 2^{-2}(1+d)^2\big((T_h T_h' - 2(1+d)^{-1}I)\gamma + T_h T_h' \tau A\big)\big((T_h T_h' - 2(1+d)^{-1}I)\gamma + T_h T_h' \tau A\big)' + 2^{-2}(1+d)^2(I - T_h T_h')(\gamma + \tau A)(\gamma + \tau A)'(I - T_h T_h')'\Big)X_*' + \delta\delta'$, $\mathfrak{R}(A)$ stands for column space of $A$ and $A^{-1}$ is an independent choice of g-inverse of $A$ and $\theta = \delta + 2^{-1}(1+d)X_*(I - T_h T_h')(\gamma + \tau A)$.

**Corollary C23:**

a) If $\lambda_* < 1$, $\hat{y}_{SRPCR}$ is superior to $\hat{y}_{SRAULE}$ in MSEM sense when the regression model is misspecified due to excluding relevant variables if and only if

$$\big((T_h T_h' - I)\gamma + T_h T_h' \tau A\big)'\Big(2^{-4}(1+d)^2(3-d)^2\tau - T_h T_h' \tau T_h' T_h)\sigma^2 + 2^{-4}\big((1+d)(3-d)\tau A - (1-d)^2\gamma\big)\big((1+d)(3-d)\tau A - (1-d)^2\gamma\big)'\Big)^{-1}\big((T_h T_h' - I)\gamma + T_h T_h' \tau A\big) \leq 1$$

where $\lambda_*$ is the largest eigenvalue of $2^4(1+d)^{-2}(3-d)^{-2}T_h T_h' \tau T_h' T_h \tau^{-1}$.

**Proof:** Consider
$$D_{(i,j)} = D(\hat{\gamma}_{SRLAUE}) - D(\hat{\gamma}_{SRPCR}) = 2^{-4}(1+d)^2(3-d)^2\sigma^2\tau - \sigma^2 T_h T_h' \tau T_h' T_h$$
$$= (2^{-4}(1+d)^2(3-d)^2\tau - T_h T_h' \tau T_h' T_h)\sigma^2$$

Since $\tau > 0$, according to Lemma A2 (see Appendix A) $D_{(i,j)} > 0$ if $\lambda_* < 1$, where $\lambda_*$ is the largest eigenvalue of $2^4(1+d)^{-2}(3-d)^{-2}T_h T_h' \tau T_h' T_h \tau^{-1}$. This completes the proof.

b) If $A \geq 0$, $\hat{y}_{SRPCR}$ is superior to $\hat{y}_{SRAULE}$ in MSEM sense when the regression model is misspecified due to excluding relevant variables if and only if $\theta \in \Re(A)$ and $\theta' A^{-1} \theta \leq 1$, where $A = X_* \Big((2^{-4}(1+d)^2(3-d)^2\tau - T_h T_h' \tau T_h' T_h)\sigma^2 + 2^{-4}((1+d)(3-d)\tau A - (1-d)^2\gamma)((1+d)(3-d)\tau A - (1-d)^2\gamma)' - ((T_h T_h' - I)\gamma + T_h T_h' \tau A)((T_h T_h' - I)\gamma + T_h T_h' \tau A)' + 2^{-4}((1+d)(3-d)I - 2^2 T_h T_h')(\gamma + \tau A)(\gamma + \tau A)'((1+d)(3-d)I - 2^2 T_h T_h')'\Big)X_*' + \delta\delta'$, $\Re(A)$ stands for column space of $A$ and $A^{-1}$ is an independent choice of g-inverse of $A$ and $\theta = \delta + 2^{-2}X_*((1+d)(3-d)I - 2^2 T_h T_h')(\gamma + \tau A)$.

**Corollary C24:**
a) If $\lambda_* < 1$, $\hat{\gamma}_{SRrk}$ is superior to $\hat{\gamma}_{SRAULE}$ in MSEM sense when the regression model is misspecified due to excluding relevant variables if and only if

$(1+k)^{-2}\big((T_h T_h' - (1+k)I)\gamma + T_h T_h' \tau A\big)' \Big((2^{-4}(1+d)^2(3-d)^2\tau - (1+k)^{-2}T_h T_h' \tau T_h' T_h)\sigma^2 +$

$2^{-4}\big((1+d)(3-d)\tau A - (1-d)^2\gamma\big)\big((1+d)(3-d)\tau A - (1-d)^2\gamma\big)'\Big)^{-1} \big((T_h T_h' - (1+k)I)\gamma +$

$T_h T_h' \tau A\big) \leq 1$

where $\lambda_*$ is the largest eigenvalue of $2^4(1+d)^{-2}(3-d)^{-2}(1+k)^{-2}T_h T_h' \tau T_h' T_h \tau^{-1}$.

**Proof:** Consider $D_{(i,j)} = D(\hat{\gamma}_{SRLAUE}) - D(\hat{\gamma}_{SRrk}) = 2^{-4}(1+d)^2(3-d)^2\sigma^2\tau - (1+k)^{-2}\sigma^2 T_h T_h' \tau T_h' T_h$
$= (2^{-4}(1+d)^2(3-d)^2\tau - (1+k)^{-2}T_h T_h' \tau T_h' T_h)\sigma^2$

Since $\tau > 0$, according to Lemma A2 (see Appendix A) $D_{(i,j)} > 0$ if $\lambda_* < 1$, where $\lambda_*$ is the largest eigenvalue of $2^4(1+d)^{-2}(3-d)^{-2}(1+k)^{-2}T_h T_h' \tau T_h' T_h \tau^{-1}$. This completes the proof.

b) If $A \geq 0$, $\hat{y}_{SRrk}$ is superior to $\hat{y}_{SRAULE}$ in MSEM sense when the regression model is misspecified due to excluding relevant variables if and only if $\theta \in \Re(A)$ and $\theta' A^{-1} \theta \leq 1$, where $A = X_* \Big((2^{-4}(1+d)^2(3-d)^2\tau - (1+k)^{-2}T_h T_h' \tau T_h' T_h)\sigma^2 + 2^{-4}((1+d)(3-d)\tau A - (1-d)^2\gamma)((1+d)(3-d)\tau A - (1-d)^2\gamma)' - (1+k)^{-2}((T_h T_h' - (1+k)I)\gamma + T_h T_h' \tau A)((T_h T_h' - (1+k)I)\gamma + T_h T_h' \tau A)' + 2^{-4}(1+k)^{-2}((1+k)(1+d)(3-d)I - 2^2 T_h T_h')(\gamma + \tau A)(\gamma + \tau A)'((1+k)(1+d)(3-d)I - 2^2 T_h T_h')'\Big)X_*' + \delta\delta'$, $\Re(A)$ stands for column space of $A$ and $A^{-1}$ is an independent choice of g-inverse of $A$ and $\theta = \delta + 2^{-2}(1+k)^{-1}X_*((1+k)(1+d)(3-d)I - 2^2 T_h T_h')(\gamma + \tau A)$.

**Corollary C25:**

a) If $\lambda_* < 1$, $\hat{\gamma}_{SRrd}$ is superior to $\hat{\gamma}_{SRAULE}$ in MSEM sense when the regression model is misspecified due to excluding relevant variables if and only if

$$2^2(1+d)^2\big((T_hT_h' - 2(1+d)^{-1}I)\gamma + T_hT_h'\tau A\big)'\big(((3-d)^2\tau - 2^2T_hT_h'\tau T_h'T_h)(1+d)^2\sigma^2 + \big((1+d)(3-d)\tau A - (1-d)^2\gamma\big)\big((1+d)(3-d)\tau A - (1-d)^2\gamma\big)'\big)^{-1}\big((T_hT_h' - 2(1+d)^{-1}I)\gamma + T_hT_h'\tau A\big) \leq 1$$

where $\lambda_*$ is the largest eigenvalue of $2^2(3-d)^{-2}T_hT_h'\tau T_h'T_h\tau^{-1}$.

**Proof:** Consider

$$D_{(i,j)} = D(\hat{\gamma}_{SRLAUE}) - D(\hat{\gamma}_{SRrd}) = 2^{-4}(1+d)^2(3-d)^2\sigma^2\tau - 2^{-2}(1+d)^2\sigma^2 T_hT_h'\tau T_h'T_h$$
$$= ((3-d)^2\tau - 2^2 T_hT_h'\tau T_h'T_h)2^{-4}(1+d)^2\sigma^2$$

Since $\tau > 0$, according to Lemma A2 (see Appendix A) $D_{(i,j)} > 0$ if $\lambda_* < 1$, where $\lambda_*$ is the largest eigenvalue of $2^2(3-d)^{-2}T_hT_h'\tau T_h'T_h\tau^{-1}$. This completes the proof.

b) If $A \geq 0$, $\hat{y}_{SRrd}$ is superior to $\hat{y}_{SRAULE}$ in MSEM sense when the regression model is misspecified due to excluding relevant variables if and only if $\theta \in \mathfrak{R}(A)$ and $\theta' A^{-1}\theta \leq 1$, where $A = X_*\bigg(((3-d)^2\tau - 2^2T_hT_h'\tau T_h'T_h)2^{-4}(1+d)^2\sigma^2 + 2^{-4}((1+d)(3-d)\tau A - (1-d)^2\gamma)((1+d)(3-d)\tau A - (1-d)^2\gamma)' - 2^{-2}(1+d)^2\big((T_hT_h' - 2(1+d)^{-1}I)\gamma + T_hT_h'\tau A\big)\big((T_hT_h' - 2(1+d)^{-1}I)\gamma + T_hT_h'\tau A\big)' + 2^{-4}(1+d)^2((3-d)I - 2T_hT_h')(\gamma + \tau A)(\gamma + \tau A)'((3-d)I - 2T_hT_h')'\bigg)X_*' + \delta\delta'$, $\mathfrak{R}(A)$ stands for column space of $A$ and $A^{-1}$ is an independent choice of g-inverse of $A$ and $\theta = \delta + 2^{-2}(1+d)X_*((3-d)I - 2T_hT_h')(\gamma + \tau A)$.

**Corollary C26:**

a) If $T_hT_h'\tau T_h'T_h$ is positive definite, $\hat{\gamma}_{SRrk}$ is superior to $\hat{\gamma}_{SRPCR}$ in MSEM sense when the regression model is misspecified due to excluding relevant variables if and only if

$$(1+k)^{-2}\big((T_hT_h' - (1+k)I)\gamma + T_hT_h'\tau A\big)'\bigg(k(2+k)(1+k)^{-2}\sigma^2 T_hT_h'\tau T_h'T_h + \big((T_hT_h' - I)\gamma + T_hT_h'\tau A\big)\big((T_hT_h' - I)\gamma + T_hT_h'\tau A\big)'\bigg)^{-1}\big((T_hT_h' - (1+k)I)\gamma + T_hT_h'\tau A\big) \leq 1$$

**Proof:** Consider $D_{(i,j)} = D(\hat{\gamma}_{SRPCR}) - D(\hat{\gamma}_{SRrk}) = \sigma^2 T_hT_h'\tau T_h'T_h - (1+k)^{-2}\sigma^2 T_hT_h'\tau T_h'T_h$
$$= \sigma^2 T_hT_h'(\tau - (1+k)^{-2}\tau)T_h'T_h$$
$$= \sigma^2 T_hT_h'k(2+k)(1+k)^{-2}\tau T_h'T_h$$
$$= k(2+k)(1+k)^{-2}\sigma^2 T_hT_h'\tau T_h'T_h$$

Since $k > 0$ and $\tau > 0$. $D_{(i,j)} > 0$ if $T_hT_h'\tau T_h'T_h$ is positive definite. This completes the proof.

b) If $A \geq 0$, $\hat{y}_{SRrk}$ is superior to $\hat{y}_{SRPCR}$ in MSEM sense when the regression model is misspecified due to excluding relevant variables if and only if $\theta \in \mathfrak{R}(A)$ and $\theta' A^{-1}\theta \leq 1$, where $A = X_*\bigg(k(2+k)(1+k)^{-2}\sigma^2 T_hT_h'\tau T_h'T_h + \big((T_hT_h' - I)\gamma + T_hT_h'\tau A\big)\big((T_hT_h' - I)\gamma + T_hT_h'\tau A\big)' - (1+k)^{-2}\big((T_hT_h' - (1+k)I)\gamma + T_hT_h'\tau A\big)\big((T_hT_h' - (1+k)I)\gamma + T_hT_h'\tau A\big)' + k^2(1+k)^{-2}T_hT_h'(\gamma + \tau A)(\gamma + \tau A)'T_h'T_h\bigg)X_*' + \delta\delta'$, $\mathfrak{R}(A)$ stands for column space of $A$ and $A^{-1}$ is an independent choice of g-inverse of $A$ and $\theta = \delta + k(1+k)^{-1}X_*T_hT_h'(\gamma + \tau A)$.

**Corollary C27:**

a) If $T_h T_h' \tau T_h' T_h$ is positive definite, $\hat{\gamma}_{SRrd}$ is superior to $\hat{\gamma}_{SRPCR}$ in MSEM sense when the regression model is misspecified due to excluding relevant variables if and only if

$2^{-2}(1+d)\big((T_h T_h' - 2(1+d)^{-1}I)\gamma + T_h T_h' \tau A\big)' \big(2^{-2}(3+d)(1-d)\sigma^2 T_h T_h' \tau T_h' T_h + \big((T_h T_h' - I)\gamma + T_h T_h' \tau A\big)\big((T_h T_h' - I)\gamma + T_h T_h' \tau A\big)'\big)^{-1}(1+d)\big((T_h T_h' - 2(1+d)^{-1}I)\gamma + T_h T_h' \tau A\big) \le 1$

**Proof:** Consider $D_{(i,j)} = D(\hat{\gamma}_{SRPCR}) - D(\hat{\gamma}_{SRrk}) = \sigma^2 T_h T_h' \tau T_h' T_h - 2^{-2}(1+d)^2 \sigma^2 T_h T_h' \tau T_h' T_h$

$= \sigma^2 T_h T_h' (\tau - 2^{-2}(1+d)^2 \tau) T_h' T_h$

$= \sigma^2 T_h T_h' 2^{-2}(3+d)(1-d)\tau T_h' T_h$

$= 2^{-2}(3+d)(1-d)\sigma^2 T_h T_h' \tau T_h' T_h$

Since $0 < d < 1$ and $\tau > 0$. $D_{(i,j)} > 0$ if $T_h T_h' \tau T_h' T_h$ is positive definite. This completes the proof.

b) If $A \ge 0$, $\hat{\gamma}_{SRrd}$ is superior to $\hat{\gamma}_{SRPCR}$ in MSEM sense when the regression model is misspecified due to excluding relevant variables if and only if $\theta \in \Re(A)$ and $\theta' A^{-1} \theta \le 1$, where $A = X_* \big(2^{-2}(3+d)(1-d)\sigma^2 T_h T_h' \tau T_h' T_h + \big((T_h T_h' - I)\gamma + T_h T_h' \tau A\big)\big((T_h T_h' - I)\gamma + T_h T_h' \tau A\big)' - 2^{-2}(1+d)^2\big((T_h T_h' - 2(1+d)^{-1}I)\gamma + T_h T_h' \tau A\big)\big((T_h T_h' - 2(1+d)^{-1}I)\gamma + T_h T_h' \tau A\big)' + 2^{-2}(1-d)^2 T_h T_h' (\gamma + \tau A)(\gamma + \tau A)' T_h' T_h\big) X_*' + \delta\delta'$, $\Re(A)$ stands for column space of $A$ and $A^{-1}$ is an independent choice of g-inverse of $A$ and $\theta = \delta + 2^{-1}(1-d)X_* T_h T_h'(\gamma + \tau A)$.

**Corollary C28:**

a) If $(k(1+d)+d-1)T_h T_h' \tau T_h' T_h$ is positive definite, $\hat{\gamma}_{SRrk}$ is superior to $\hat{\gamma}_{SRrd}$ in MSEM sense when the regression model is misspecified due to excluding relevant variables if and only if

$(1+k)^{-2}\big((T_h T_h' - (1+k)I)\gamma + T_h T_h' \tau A\big)' \big(2^{-2}(1+k)^{-2}(k(1+d)+d+3)(k(1+d)+d-1)\sigma^2 T_h T_h' \tau T_h' T_h + 2^{-2}(1+d)^2\big((T_h T_h' - 2(1+d)^{-1}I)\gamma + T_h T_h' \tau A\big)\big((T_h T_h' - 2(1+d)^{-1}I)\gamma + T_h T_h' \tau A\big)'\big)^{-1} \big((T_h T_h' - (1+k)I)\gamma + T_h T_h' \tau A\big) \le 1$

**Proof:** Consider $D_{(i,j)} = D(\hat{\gamma}_{SRrd}) - D(\hat{\gamma}_{SRrk}) = 2^{-2}(1+d)^2 \sigma^2 T_h T_h' \tau T_h' T_h - (1+k)^{-2} \sigma^2 T_h T_h' \tau T_h' T_h$

$= \sigma^2 T_h T_h' (2^{-2}(1+d)^2 - (1+k)^{-2})\tau T_h' T_h$

$= \sigma^2 T_h T_h' 2^{-2}(k(1+d)+d+3)(k(1+d)+d-1)(1+k)^{-2}\tau T_h' T_h$

$= 2^{-2}(1+k)^{-2}\sigma^2 (k(1+d)+d+3)(k(1+d)+d-1)T_h T_h' \tau T_h' T_h$

Since $0 < d < 1$, $k > 0$ and $\tau > 0$. $D_{(i,j)} > 0$ if $(k(1+d)+d-1)T_h T_h' \tau T_h' T_h$ is positive definite. This completes the proof.

b) If $A \ge 0$, $\hat{\gamma}_{SRrk}$ is superior to $\hat{\gamma}_{SRrd}$ in MSEM sense when the regression model is misspecified due to excluding relevant variables if and only if $\theta \in \Re(A)$ and $\theta' A^{-1} \theta \le 1$, where $A = X_* \big(2^{-2}(1+k)^{-2}\sigma^2(k(1+d)+d+3)(k(1+d)+d-1)T_h T_h' \tau T_h' T_h + 2^{-2}(1+d)^2\big((T_h T_h' - 2(1+d)^{-1}I)\gamma + T_h T_h' \tau A\big)\big((T_h T_h' - 2(1+d)^{-1}I)\gamma + T_h T_h' \tau A\big)' - (1+k)^{-2}\big((T_h T_h' - (1+k)I)\gamma + T_h T_h' \tau A\big)\big((T_h T_h' - (1+k)I)\gamma + T_h T_h' \tau A\big)' + 2^{-2}(1+k)^{-2}\big((1+d)(1+k) - 2\big)^2 T_h T_h'(\gamma + \tau A)(\gamma + \tau A)' T_h' T_h\big) X_*' + \delta\delta'$, $\Re(A)$ stands for column space of $A$ and $A^{-1}$ is an independent choice of g-inverse of $A$ and $\theta = \delta + 2^{-1}(1+k)^{-1}((1+d)(1+k) - 2)X_* T_h T_h'(\gamma + \tau A)$.